\input amstex
\documentstyle{amsppt}

\magnification=\magstep1
\hsize=5.2in
\vsize=6.7in

\catcode`\@=11
\loadmathfont{rsfs}
\def\mycal{\mathfont@\rsfs}
\csname rsfs \endcsname

\topmatter
\title  SOME COMPUTATIONS OF $1$-COHOMOLOGY GROUPS \\
AND CONSTRUCTION OF NON ORBIT EQUIVALENT ACTIONS
\endtitle
\author SORIN POPA \endauthor

\rightheadtext{Computations of 1-cohomology groups}

\affil University of California, Los Angeles\endaffil

\address Math.Dept., UCLA, LA, CA 90095-155505\endaddress
\email popa@math.ucla.edu\endemail

\thanks Supported in part by NSF Grant 0100883.\endthanks

\abstract For each group $G$ having an infinite normal subgroup
with the relative property (T) (e.g.
$G=H \times K$, with $H$ infinite with property (T) and $K$
arbitrary) and each countable abelian group $\Lambda$ we construct
free ergodic measure-preserving actions $\sigma_\Lambda$ of $G$ on
the probability space such that the $1$'st cohomology group of
$\sigma_\Lambda$, $\text{\rm H}^1(\sigma_\Lambda,G)$, is equal to
Char$(G) \times \Lambda$. We deduce that $G$ has uncountably  many
non stably orbit equivalent actions. We also calculate
1-cohomology groups and show existence of ``many'' non stably
orbit equivalent actions for free products of groups as above.
\endabstract

\endtopmatter

\document

\heading 0. Introduction
\endheading

Let $G$ be a countable discrete group and $\sigma: G \rightarrow
{\text{\rm Aut}}(X,\mu)$ a free measure preserving (m.p.) action
of $G$ on the probability space $(X, \mu)$, which we also view as
an integral preserving action of $G$ on the abelian von Neumann
algebra $A=L^\infty(X,\mu)$. A 1-{\it cocycle} for $(\sigma, G)$
is a map $w: G \rightarrow \Cal U(A)$, satisfying
$w_g\sigma_g(w_h)=w_{gh}, \forall g,h \in G$, where $\Cal
U(A)=\{u\in A\mid uu^*=1\}$ denotes the group of unitary elements
in $A$. The set of 1-cocycles for $\sigma$ is denoted $\text{\rm
Z}^1(\sigma,G)$ and is endowed with the Polish group structure
given by point multiplication and pointwise convergence in the
norm $\|\cdot\|_2$. The 1-{\it cohomology group} of $\sigma$,
$\text{\rm H}^1(\sigma,G)$, is the quotient of $\text{\rm
Z}^1(\sigma,G)$ by the subgroup of coboundaries $\text{\rm
B}^1(\sigma,G)=\{\sigma_g(u)u^* \mid u\in \Cal U(A)\}$.

The group $\text{\rm H}^1(\sigma,G)$ was first mentioned by I.M.
Singer  ([Si55]), related to his study of automorphisms of group
measure space von Neumann algebras. J. Feldman and C.C. Moore
extended the definition to countable, measurable equivalence
relations and pointed out that $\text{\rm H}^1(\sigma,G)$ depends
only on the orbit equivalence (OE) class of $(\sigma,G)$ ([FM77]),
thus being an OE invariant for actions. K. Schmidt showed in
([S80], [S81]) that $\text{\rm H}^1(\sigma,G)$ is Polish (i.e.
$\text{\rm B}^1(\sigma, G)$ closed in $\text{\rm Z}^1(\sigma,G)$)
if and only if $\sigma$ has no non-trivial asymptotically
invariant sequences and noticed that Bernoulli shift actions of
non-amenable groups are always strongly ergodic (so their H$^1$
group is Polish). On the other hand, by ([D63], [OW80], [CFW81])
all free ergodic m.p. actions of infinite amenable groups are OE
and non-strongly ergodic, thus having all the same (``wild'')
H$^1$-group. Results of A. Connes and B. Weiss in ([CW80]) and
([S81]) show that $\text{\rm H}^1(\sigma,G)$ is countable discrete
$\forall \sigma$ if and only if $G$ has the property (T) of
Kazhdan.

C.C. More produced the first examples of free ergodic m.p. actions
of infinite groups with trivial 1-cohomology  ([M82] page 220;
note that the groups in these examples have the property (T) of
Kazhdan). Then in ([Ge87]) Gefter showed that if a Kazhdan group
$G$ can be densely embedded into a compact simply connected
semi-simple Lie group $\Cal G$ and $K \subset \Cal G$ is a closed
subgroup then the action of $G$ by left translation on $\Cal G/K$
has H$^1$-group equal to Char$(G \times K)$ (the character group
of $G \times K$). But these initial calculations were not followed
up upon and, in fact, after an intense activity during 1977-1987
([FM77], [OW80], [S80], [S81], [CFW81], [CW80], [Z84], [Ge87],
[GeGo87], [P86], [JS87]), the whole area of orbit equivalence
ergodic theory went through more than a decade of relative
neglect.

But even after the spectacular revival of this subject in recent
years ([Fu99], [G00], [G02], [Hj02], [P02], [MoSh02], [GP03],
[P04]), 1-cohomology wasn't really exploited as a tool to
distinguish between orbit inequivalent actions of groups. And this
despite a new calculation of H$^1$-groups was obtained
in ([PSa03]), this time for Bernoulli shift actions $\sigma$ of
ANY property (T) groups, more generally for groups $G$ having
infinite normal subgroups with the relative property (T) of
Kazhdan-Margulis (called {\it weakly rigid} in [P01], [P03],
[PSa03]). Thus, it is shown in ([PSa03]) that for such
$(\sigma,G)$ one has  $\text{\rm H}^1(\sigma,G)=\text{\rm
Char}(G)$.

In this paper we consider an even larger class of groups, denoted
$w\Cal T$, generalizing the weakly rigid groups of ([P01],
[PSa03]), and for each $G \in w\Cal T$ calculate H$^1(\sigma',G)$
for a large family of quotients $\sigma': G \rightarrow \text{\rm
Aut}(X',\mu')$ of the Bernoulli shift actions $\sigma$ of $G$ on
$(X, \mu) = \Pi_{g\in G} (\Bbb T, \lambda)_g$. Thus, our main
result shows that given any countable discrete abelian group
$\Lambda$ there exists a free ergodic action $\sigma_\Lambda$ of
$G$ implemented by the restriction of $\sigma$ to an appropriate
$\sigma$-invariant subalgebra of $L^\infty(X,\mu)$, such that
H$^1(\sigma_\Lambda,G) = \text{\rm Char}(G) \times \Lambda$ as
topological groups. We also calculate the 1-cohomology for similar
quotients of Bernoulli shift actions of free products of groups in
$w\Cal T$. We deduce that each $G \in w\Cal T$, or $G$ a free
product of infinite Kazhdan groups, has a continuous family of
free ergodic m.p. actions with mutually non-isomorphic
H$^1$-groups, and which are thus OE inequivalent.

These results, together with prior ones in ([Ge87], [PSa03]),
establish 1-cohomology as an effective OE invariant, adding to the
existing pool of methods used to differentiate orbit inequivalent
actions ([CoZ89], [Fu99], [G00], [G02], [MoSh02]). Rather than
depending on the group only, like the cost or $\ell^2$-Betti
numbers in ([G00], [G02]), the H$^1$-invariant depends also on the
action, proving particularly useful in distinguishing large
classes of orbit inequivalent actions of the same group.

Before stating the results in more details, let us define
the class $w\Cal T$ more precisely:
It consists of all countable groups $G$ which
contain an infinite subgroup $H \subset G$ with the relative
property (T) of Kazhdan-Margulis (see [Ma82], [dHV89]) such that
$H$ is {\it wq-normal} in $G$, i.e. given any intermediate
subgroup $H\subset K \varsubsetneq G$ there exists $g \in
G\setminus K$ with $gKg^{-1} \cap K$ infinite (see 2.3 for other
equivalent characterizations of this property). For instance, if
there exist finitely many subgroups  $H=H_0 \subset H_1 \subset
... \subset H_n =G$ with all consecutive inclusions $H_j \subset
H_{j + 1}$ normal, then $H \subset G$ is wq-normal. In particular,
weakly rigid groups are in the class $w\Cal T$.

\proclaim{Theorem 1} Let $G\in w\Cal T$. Let $\sigma$ be a
Bernoulli shift action of $G$ on the probability space $(X, \mu)=
\Pi_g(X_0, \mu_0)_g$ and $\beta$ a free action of a group $\Gamma$
on $(X, \mu)$ that commutes with $\sigma$ and such that the
restriction $\sigma^\Gamma$ of $\sigma$ to the fixed point algebra
$\{a\in L^\infty(X, \mu) \mid \beta_h(a)=a, \forall h\in \Gamma
\}$ is still a free action of $G$. Denote $\text{\rm
Char}_\beta(\Gamma)$ the group of characters $\gamma$ on $\Gamma$
for which there exist unitary elements $u \in L^\infty(X,\mu)$
with $\beta_h(u)=\gamma(h)u, \forall h\in \Gamma$. Then
$\text{\rm H}^1(\sigma^\Gamma, G) = \text{\rm Char}(G) \times \text{\rm
Char}_\beta(\Gamma)$ as topological groups.
\endproclaim

Since any countable abelian group $\Lambda$ can be realized as
$\text{\rm Char}_\beta(\Gamma)$, for some appropriate action
$\beta$ of a group $\Gamma$ commuting with the Bernoulli shift
$\sigma$, and noticing that $\text{\rm H}^1(\sigma,G)$ are even
invariant to stable orbit equivalence, we deduce:

\proclaim{Corollary 2} Let $G\in w\Cal T$. Given any countable
discrete abelian group $\Lambda$ there exists a free ergodic m.p.
action $\sigma_\Lambda$ of $G$ on the standard non-atomic
probability space such that $\text{\rm H}^1(\sigma_\Lambda,G)=\text{\rm
Char}(G) \times \Lambda$. Moreover, $\sigma_\Lambda$ can be taken
``quotients'' of $G$-Bernoulli shifts. Thus, any $G\in w\Cal T$
has a continuous family of mutually non-stably orbit equivalent
free ergodic m.p. actions on the probability space.
\endproclaim

Examples of groups in the class $w\Cal T$ covered by the above
results are the infinite property (T) groups, the groups $\Bbb Z^2
\rtimes \Gamma$, for $\Gamma \subset SL(2, \Bbb Z)$ non-amenable
(cf [K67], [Ma82], [B91]) and the groups $\Bbb Z^N \rtimes \Gamma$
for suitable actions of arithmetic lattices $\Gamma$ in $SU(n, 1)$
or $SO(n,1), n \geq 2$ (cf [V04]). Note that if $G \in w\Cal T$
and $K$ is a group acting on $G$ by automorphisms then $G \rtimes
K \in w\Cal T$. Also, if $G\in w\Cal T$ and $K$ is an arbitrary
group then $G \times K \in w\Cal T$. In particular, any product
between an infinite property (T) group and an arbitrary group is
in the class $w\Cal T$. Thus, Corollary 2 covers a recent result
of G. Hjorth ([Hj02]), showing that infinite property (T) groups
have uncountably many orbit inequivalent actions. Moreover, rather
than an existence result, Corollary 2 provides a concrete list of
uncountably many
inequivalent actions
(indexed by the virtual isomorphism classes
of all countable, discrete,
abelian groups, see 2.13).

Note however that if $G = H * H'$ with $H$ infinite with property
(T) and $H'$ non-trivial, then $(G,H)$ does have the relative
property (T) but $H$ is not wq-normal in $G$. One can in fact show
that $G$ is not in the class $w\Cal T$. Yet we can still calculate
in this case the 1-cohomology for the quotients of $G$-Bernoulli
shifts $\sigma_\Lambda$ considered in Corollary 2. While
$\text{\rm H}^1(\sigma_\Lambda,G)$ are ``huge'' (non locally
compact) in this case, if we denote by $\tilde{\text{\rm
H}}^1(\sigma_\Lambda,G)$ the quotient of $\text{\rm
H}^1(\sigma_\Lambda,G)$ by the connected component of $1$ then we
get:

\proclaim{Theorem 3} Let $\{G_n\}_{n \geq 0}$ be a sequence of
countable groups such that each $G_n$ is either amenable or
belongs to the class $w\Cal T$ and denote $G=*_{n\geq 0} G_n$.
Assume the set $J$ of indices
$j\geq 0$ for which $G_j\in w\Cal T$ is non-empty and that
$G_j$ has totally disconnected character group, $\forall j\in J$.
If $\Lambda$ is a countable abelian group, then
$\tilde{\text{\rm H}}^1(\sigma_\Lambda, G) \simeq \Pi_{j \in J} \text{\rm
Char}(G_j) \times \Lambda^{|J|}$ as Polish groups.
\endproclaim

Since property (T) groups have  finite
(thus totally disconnected) character group, from the above
theorem we get:

\proclaim{Corollary 4} Let $H_1, H_2, ..., H_k$ be infinite
property $(\text{\rm T})$ groups and $0 \leq n \leq \infty$. The
free product group $H_1 * H_2 * ... * H_k * \Bbb F_n$ has
continuously many non stably orbit equivalent free ergodic m.p.
actions.
\endproclaim

The use of von Neumann algebras framework and non-commutative
analysis tools is crucial for the approach in this paper. Thus,
the construction used in Theorem 1, as well as its proof, become
quite natural in von Neumann algebra context, where similar ideas
have been used in ([P01]), to compute the 1-cohomology and
fundamental group for non-commutative (Connes-St\o rmer) Bernoulli
shift actions of weakly rigid groups on the hyperfinite II$_1$
factor $R$, and in ([C75b]), to compute the approximately inner,
centrally free part $\Cal X (M)$ of the outer automorphism group
of a II$_1$ factor $M$.

The paper is organized as follows: In Section 1 we present some
basic facts on 1-cohomology for actions, including a detailed
discussion of the similar concept for full groupoids and
equivalence relations. Also, we revisit the results on
1-cohomology in ([FM77], [S80], [S81]). In Section 2 we prove
Theorem 1 and its consequences. In Section 3 we consider actions
of free product groups and prove Theorem 3.

This work was done while I was visiting the Laboratoires
d'Alg\'ebres  d'Op\'erateurs at the Instituts de Math\'emathiques
of the Universities Luminy and Paris 7 during the Summer of 2004.
I am grateful to the CNRS and the members of these Labs for their
support and kind hospitality. It is a pleasure for me to thank
Bachir Bekka, Etienne Ghys, Alekos Kechris and Stefaan Vaes for
generous comments and useful discussions. I am particularly
grateful to Sergey Gefter and Antony Wasserman for pointing out to
me the calculations of 1-cohomology groups in ([M82], [Ge87]).

\heading 1. $1$-cohomology for actions and equivalence relations
\endheading

We recall here the definition and basic properties of the
1-cohomology groups for actions and equivalence relations, using
the framework of von Neumann algebras. We revisit this way the
results in ([S80], [S81], [FM77]) and prove the invariance of
1-cohomology groups to stable orbit equivalence. The von Neumann
algebra setting leads us to adopt H. Dye's initial point of view
([D63]) of regarding equivalence relations as full groupoids and
to use I.M. Singer's observation ([Si55]) that the group of
1-cocycles of an action is naturally isomorphic to the group of
automorphisms of the associated group measure space von Neumann
algebra that leave the Cartan subalgebra pointwise fixed.

\vskip .05in
\noindent {\it 1.1. $1$-cohomology for actions}. Let $\sigma: G
\rightarrow {\text{\rm Aut}}(X, \mu)$ be a free measure preserving
(abbreviated m.p.) action of the (at most) countable discrete
group $G$ on the standard probability space $(X, \mu)$ and still
denote by $\sigma$ the action it implements on $A = L^\infty(X,
\mu)$. Denote $\Cal U(A)=\{u\in A \mid uu^*=1\}$ the group of
unitary elements of $A$. A function $w : G \rightarrow \Cal U(A)$
satisfying $w_g\sigma_g(w_h)=w_{gh}, \forall g,h \in G$, is called
a $1$-{\it cocycle} for $\sigma$. Note that a scalar valued
function $w: G \rightarrow \Cal U(A)$ is a $1$-cocycle iff $w \in
{\text{\rm Char}}(G)$.

Two $1$-cocycles $w, w'$ are {\it cohomologous}, $w \sim_c w'$, if
there exists $u \in \Cal U(A)$ such that $w_g' = u^* w_g
\sigma_g(u), \forall g\in G$. A 1-cocycle $w$ is {\it coboundary}
if $w \sim_c {\bold 1}$, where ${\bold 1}_g =1, \forall g$.

Denote by $\text{\rm Z}^1(\sigma, G)$ (or simply $\text{\rm Z}^1(\sigma)$,
when there is no risk
of confusion) the set of
$1$-cocycles for $\sigma$, endowed with the structure of a
topological (commutative) group given by point multiplication and
pointwise convergence in norm $\|\cdot\|_2$. Denote by
$\text{\rm B}^1(\sigma, G)=\text{\rm B}^1(\sigma)\subset
\text{\rm Z}^1(\sigma)$ the subgroup of
coboundaries and by $\text{\rm H}^1(\sigma, G)=
\text{\rm H}^1(\sigma)$ the quotient
group $\text{\rm Z}^1(\sigma)/\text{\rm B}^1(\sigma)
= \text{\rm Z}^1(\sigma)/\sim_c$, called the
$1$'{\it st cohomology group} of $\sigma$. Note that $\text{\rm
Char}(G)$ with its usual topology can be viewed as a compact
subgroup of $\text{\rm Z}^1(\sigma)$ and its image in
$\text{\rm H}^1(\sigma)$ is a
compact subgroup.  If in addition $\sigma$ is weakly mixing, then
the image of Char$(G)$ in $\text{\rm H}^1(\sigma)$ is faithful (see
2.4.1$^\circ$).

The groups $\text{\rm B}^1(\sigma),\text{\rm Z}^1(\sigma),
\text{\rm H}^1(\sigma)$ were first considered in ([Si55]). As
noticed in ([Si55]), they can be identified with certain groups of
automorphisms of the finite von Neumann algebra $M = A
\rtimes_\sigma G$, as explained below. Note that there exists a
unique normal faithful trace $\tau$ on $M$ that extends the
integral $\int \cdot {\text{\rm d}}\mu$ on $A$ and that $M$ is a
factor iff $\sigma$ is ergodic. For $x\in M$ we denote
$\|x\|_2=\tau(x^*x)^{1/2}$.

\vskip .05in \noindent {\it 1.2. Some related groups of
automorphisms}. Let Aut$_0(M;A)$ denote the  group of
automorphisms of $M$ that leave all elements of $A$ fixed, endowed
with the topology of pointwise convergence in norm $\|\cdot \|_2$
(the topology it inherits from Aut$(M,\tau)$). If $\theta \in
{\text{\rm Aut}}_0(M;A)$ then $w^\theta_g = \theta(u_g)u_g^*$, $g
\in G$, is a 1-cocycle, where $\{u_g\}_g\subset M$ denote the
canonical unitaries implementing the action $\sigma$. Conversely,
if $w \in \text{\rm Z}^1(\sigma)$ then $\theta^w(au_g)=aw_gu_g, a \in A, g\in
G$, defines an automorphism of $M$ that fixes $A$. Clearly $\theta
\mapsto w^\theta$, $w \mapsto \theta^w$ are group morphisms and
are inverse one another, thus identifying $\text{\rm Z}^1(\sigma)$ with
Aut$_0(M;A)$ as topological groups, with $\text{\rm B}^1(\sigma)$
corresponding to the inner automorphism group Int$_0(M;A) =
\{\text{\rm Ad}(u) \mid u \in \Cal U(A)\}$. Thus, $\text{\rm H}^1(\sigma)$ is
naturally isomorphic to Out$_0(M;A) \overset \text{\rm def} \to =
{\text{\rm Aut}}_0(M;A)/{\text{\rm Int}}_0(M;A)$.

The groups ${\text{\rm Aut}}_0(M;A)$, ${\text{\rm Int}}_0(M;A)$,
Out$_0(M;A)$ make actually sense for any inclusion $A \subset M$
consisting of a II$_1$ factor $M$ with a {\it Cartan subalgebra}
$A$, i.e. a maximal abelian $*$-subalgebra of $M$ with {\it
normalizer} $\Cal N_M(A)\overset \text{\rm def} \to = \{ u \in
\Cal U(M) \mid uAu^* = A \}$ generating $M$. In order to interpret
Out$_0(M;A)$ as 1-cohomology group in this more general case,
we'll recall from ([D63], [FM77]) two alternative, equivalent ways
of viewing Cartan subalgebra inclusions $A \subset M$.

\vskip .05in \noindent {\it 1.3. Full groupoids and equivalence
relations}. With $A \subset M$ as above, let $\Cal G\Cal N_M(A) =
\{v\in M \mid vv^*, v^*v \in \Cal P(M), vAv^* = Avv^*\}$, where
$\Cal P(M)$ denotes the idempotents (or projections) in $A$.
Identify $A$ with $L^\infty(X, \mu)$, for some probability space
$(X, \mu)$, with $\mu$ corresponding to $\tau_{|A}$, where $\tau$
is the trace on $M$. We denote by $\Cal G_{A \subset M}$ the set
of all local isomorphisms $\phi_v = {\text{\rm Ad}}(v)$, $v \in
\Cal G\Cal N_M(A)$, defined modulo sets of measure zero. We endow
$\Cal G_{A \subset M}$ with the natural groupoid structure given
by composition, and call it the {\it full groupoid} associated to
$A \subset M$. Since $\{v_n\}_n \subset \Cal G\Cal N_M(A)$ with
$\{v_nv_n^*\}_n$, resp. $\{v_n^*v_n\}_n$, mutually orthogonal
implies $\Sigma_n v_n \in \Cal G\Cal N_M(A)$, it follows that
$\Cal G=\Cal G_{A \subset M}$ satisfies the axiom:

\vskip .05in \noindent $(1.3.1)$. Let $R, L \subset X$ be
measurable subsets with $\mu(R)=\mu(L)$ and $\phi : R \simeq L$ a
measurable, measure preserving isomorphism. Then $\phi \in \Cal G$
iff there exists a countable partition of $R$ with measurable
subsets $\{R_n\}_n$ such that $\phi_{|R_n}\in \Cal G, \forall n$.
In particular, the set of {\it units} $\Cal G_0=\{\phi^{-1} \phi
\mid \phi \in \Cal G\}$ of $\Cal G$ is equal to $\{ id_{Y} \mid Y
\subset X$ measurable $\}$. \vskip .05in

Note that the factoriality of $M$ amounts to the ergodicity of the
action of $\Cal G$ on $L^\infty(X, \mu)$ and that $M$ is separable
in the norm $\|\cdot \|_2$ iff $\Cal G$ is countably generated as
a groupoid satisfying $(1.3.1)$. If $M=A \rtimes_\sigma G$ for
some free m.p. action $\sigma$ of a group $G$, then we denote
$\Cal G_{A \subset M}$ by $\Cal G_\sigma$. Note that if $\phi : R
\simeq L$ is a m.p.  isomorphism, for some measurable
subsets $R, L \subset X$ with $\mu(R)=\mu(L)$,
then $\phi \in \Cal G_\sigma$ iff there
exist $g_n \in G$ and a partition of $R$ with measurable subsets
$\{R_n\}_n$ such that $\phi_{|R_n}=\sigma(g_n)_{|R_n}, \forall n$.

A groupoid $\Cal G$ of m.p. local
isomorphisms of the probability space $(X, \mu)$ satisfying
$(1.3.1)$ is called an {\it abstract} {\it full groupoid}.

If an abstract groupoid $\Cal G$ acting on $(X, \mu)$ is given
then let $\Bbb C \Cal G$ denote the groupoid algebra of formal
finite linear combinations $\Sigma_\phi c_\phi \phi$. Let
$\tau(\phi)$ denote the measure of the largest set on which $\phi$
acts as the identity and extend it by linearity to $\Bbb C \Cal
G$. Then define a sesquilinear form on $\Bbb C \Cal G$ by $\langle
x, y \rangle=\tau(y^*x)$ and denote by $L^2(\Cal G)$ the Hilbert
space obtained by completing $\Bbb C \Cal G/I_\tau$ in the norm
$\|x\|_2=\tau(x^*x)^{1/2}$, where $I_\tau=\{x \mid \langle x,x
\rangle =0\}$. Each $\phi\in \Cal G$ acts on $L^2(\Cal G)$ as a
left multiplication operator $u_\phi$. Denote by $L(\Cal G)$ the
von Neumann algebra generated by the operators $\{u_\phi, \phi \in
\Cal G\}$ and by $L(\Cal G_0)\simeq L^\infty(X, \mu)$ the von
Neumann subalgebra generated by the units $\Cal G_0$.

It is easy to check that $L(\Cal G)$ is a finite von Neumann
algebra with Cartan subalgebra $L(\Cal G_0)=L^\infty(X, \mu)$ and
faithful normal trace $\tau$ extending the integral on
$L^\infty(X, \mu)$ and satisfying $\tau(u_\phi)=\tau(\phi)$, and
with $L^2(L(\Cal G))=L^2(\Cal G)$ the standard representation of
$(L(\Cal G), \tau)$. Moreover, if $A = L(\Cal G_0)$, $M = L(\Cal
G)$ then $\Cal G\Cal N_M(A) = \{au_\phi \mid \phi \in \Cal G, a
\in \Cal I(A)\}$, where $\Cal I(A)$ denotes the set of partial
isometries in $A$. Thus, the measurable groupoid $\Cal G_{A
\subset M}$ associated to the Cartan subalgebra inclusion $L(\Cal
G_0) \subset L(G)$ can be naturally identified with $\Cal G$.
Also, note that $L(\Cal G)$ is a factor iff $\Cal G$ is ergodic,
in which case either $L(\Cal G)\simeq M_{n\times n}(\Bbb C)$ (when
$(X,\mu)$ is the $n$-points probability space) or $L(\Cal G)$ is a
II$_1$ factor (when $(X, \mu)$ has no atoms, equivalently when
$\Cal G$ has infinitely many elements).

If the full groupoid $\Cal G$ is generated by a countable set of
local isomorphisms $\{\phi_n\}_n \subset \Cal G$ and one considers
a standard Borel structure on $X$ with $\sigma$-field $\Cal X$
then $\phi_n$ can be taken Borel. If one denotes $\Cal R=\Cal
R_{\Cal G}$ the equivalence relation implemented by the orbits of
$\phi \in \Cal G$ then each class of equivalence in $\Cal R$ is countable and
$\Cal R$ lies in the $\sigma$-field $\Cal X \times \Cal X$.
Moreover, all $\phi \in \Cal G$ can be recuperated from $\Cal R$
as graphs of local isomorphisms that lie in $\Cal R \cap \Cal X
\times \Cal X$. Such $\Cal R$ is called a {\it countable measure preserving}
(m.p.) {\it standard equivalence relation}. The m.p. standard
equivalence relation $\Cal R_{A \subset M}$ associated to a Cartan
subalgebra inclusion $A \subset M$ is the equivalence relation
implemented by the orbits of $\Cal G_{A \subset M}$. In the case
$\Cal G$ is given by an action $\sigma$ of a countable group $G$,
the orbits of $\Cal G_\sigma$ coincide with the orbits of $\sigma$
and one denotes the corresponding equivalence relation by $\Cal
R_\sigma$.

An isomorphism between two full groupoids (resp. m.p. equivalence
relations) is an isomorphism of the corresponding probability
spaces that takes one groupoid (resp. m.p. equivalence relation)
onto the other. Such an isomorphism clearly agrees with the
correspondence between groupoids and equivalence relations
described above. Two Cartan subalgebra inclusions $(A_1 \subset
M_1, \tau_1)$, $(A_2 \subset M_2, \tau_2)$  are isomorphic if
there exists $\theta : (M_1, \tau_1) \simeq (M_2, \tau_2)$ such
that $\theta(A_1)=A_2$. Note that if this is the case then $\Cal
G_{A_1\subset M_1} \simeq \Cal G_{A_2 \subset M_2}$, $\Cal
R_{A_1\subset M_1} \simeq \Cal R_{A_2 \subset M_2}$. Conversely,
if $\Cal G_1 \simeq \Cal G_2$ then $(L(G_{1,0}) \subset L(G_1))
\simeq (L(G_{2,0}) \subset L(G_2))$. In particular, two free
ergodic m.p. actions $\sigma_i : G_i \rightarrow {\text{\rm
Aut}}(X_i, \mu_i)$ are orbit equivalent iff $(A_1 \subset
A_1\rtimes_{\sigma_1} G_1) \simeq (A_2 \subset
A_2\rtimes_{\sigma_2} G_2)$.

\vskip .05in \noindent {\it 1.4. Amplifications and stable orbit
equivalence}. If $M$ is a II$_1$ factor and $ t > 0$ then for any
$n \geq m \geq t$ and any projections $p\in M_{n\times n}(M)$,
$q\in M_{m\times m}(M)$ of (normalized) trace $\tau(p)=t/n$,
$\tau(q)=t/m$, one has $pM_{n\times n}(M)p \simeq qM_{m\times
m}(M)q$. Indeed, because if we regard $M_{m\times m}(M)$ as a
``corner'' of $M_{n \times n}(M)$ then $p,q$ have the same trace
in $M_{n \times n}(M)$, so they are conjugate by a unitary $U$ in
$M_{n\times n}(M)$, which implements an isomorphism between
$pM_{n\times n}(M)p$ and $qM_{m\times m}(M)q$. One denotes by
$M^t$ this common (up to isomorphism)  II$_1$ factor and callit
the {\it amplification} of $M$ by $t$.

Similarly, if $A \subset M$ is a Cartan subalgebra of the II$_1$
factor $M$ then $(A \subset M)^t = (A^t \subset M^t)$ denotes the
(isomorphism class of the) Cartan subalgebra inclusion $p(A
\otimes D_n \subset M \otimes M_{n \times n}(\Bbb C))p$ where $n
\geq t$, $D_n$ is the diagonal subalgebra of $M_{n \times n}(\Bbb
C)$ and $p \in A \otimes D_n$ is a projection of trace $\tau(p) =
t/n$. In this case, the fact that the isomorphism class of $(A
\subset M)^t$ doesn't depend on the choice of $n, p$ follows from
a lemma of H. Dye ([D63]), showing that if $M_0$ is a II$_1$
factor and $A_0 \subset M_0$ is a Cartan subalgebra, then two
projections $p,q \in A_0$ having the same trace are conjugate by a
unitary element in the normalizer of $A_0$ in $M_0$.

If $\Cal G$ is an
ergodic full groupoid on the non-atomic probability space then
$\Cal G^t$ is the full groupoid obtained by restricting the full
groupoid generated by $\Cal G \times \Cal D_n$ to a subset of
measure $t/n$, where $\Cal D_n$ is the groupoid of permutations of
the $n$-points probability space with the counting measure. If
$\Cal R$ is an ergodic m.p. standard equivalence relation then
$\Cal R^t$ is defined in a similar way. Again, $\Cal G^t$, $\Cal R^t$
are defined only up to isomorphism.

$(A \subset M)^t$ (resp. $\Cal G^t$, $\Cal R^t$) is called the
$t$-{\it amplification} of $A \subset M$ (resp. of $\Cal G$, $\Cal
R$). We clearly have $\Cal G_{(A \subset M)^t}=\Cal G_{(A \subset
M)}^t$, $\Cal R_{(A \subset M)^t}=\Cal R_{(A \subset M)}^t$ and if
$\Cal G$, $\Cal R$ correspond with one another then so do $\Cal
R^t, \Cal G^t$, $\forall t$. Note that $((A \subset M)^t)^s= (A
\subset M)^{st}$, $(\Cal G^t)^s=\Cal G^{ts}$, $(\Cal R^t)^s = \Cal
R^{ts}$, $\forall t, s >0$.

Two ergodic full groupoids $\Cal G_i, i=1,2$ (resp. ergodic
equivalence relations $\Cal R_i, i=1,2$) are {\it stably orbit
equivalent} if $\Cal G_1 \simeq \Cal G_2^t$ (resp. $\Cal R_1
\simeq \Cal R_2^t$), for some $t > 0$. Two free ergodic m.p.
actions $(\sigma_i,G_i), i=1,2$ are stably orbit equivalent if
$\Cal G_{\sigma_1} \simeq \Cal G_{\sigma_2}^t$ for some $t$. Note
that this is equivalent to the existence of subsets of positive
measure $Y_i\subset X_i$ and of an isomorphism $\Psi: (Y_1,
\mu_1/\mu_1(Y_1)) \simeq (Y_2, \mu_2/\mu_2(Y_2))$ such that
$\Psi(\sigma_1(G_1)x \cap Y_1)=\sigma_2(G_2)\Psi(x)\cap Y_2$, a.e.
in $x\in Y_1$.

\vskip .05in \noindent {\it 1.5. $1$-cohomology for full
groupoids}. Let $\Cal G$ be a full groupoid acting on the
probability space $(X, \mu)$ and denote $A = L^\infty(X, \mu)$. A
1-{\it cocycle} for $\Cal G$ is a map $w : \Cal G \rightarrow \Cal
I(A)$ satisfying the relation $w_\phi \phi(w_\psi)= w_{\phi
\psi}$, $\forall \phi, \psi \in \Cal G$. In particular, this implies that
the support of $w_\phi$, $w_\phi w_\phi^*$, is equal to the range
$r(\phi)$ of $\phi$. Thus, $w_{id_Y}=\chi_Y, \forall Y \subset X$
measurable.

We denote by $\text{\rm Z}^1(\Cal G)$ the set of all $1$-cocycles and endow
it with the (commutative) semigroup structure given by point
multiplication. We denote by ${\bold 1}$ the 1-cocycle given by
${\bold 1}_\phi = r(\phi), \forall \phi \in \Cal G$. If we let
$(w^{-1})_\phi = {w_\phi}^*$ then we clearly have $w w^{-1}=
{\bold 1}$ and ${\bold 1} w = w$, $\forall w \in \text{\rm Z}^1(\Cal G)$.
Thus, together also with the topology given by pointwise norm
$\|\cdot\|_2$-convergence, $\text{\rm Z}^1(\Cal G)$ is a commutative Polish
group.

Two 1-cocycles $w_1, w_2$ are {\it cohomologous}, $w_1 \sim_c
w_2$, if there exists $u \in \Cal U(A)$ such that $w_2(\phi) = u^*
w_2(\phi) \phi(u)$, $\forall \phi \in \Cal G$. A 1-cocycle $w$
cohomologous to ${\bold 1}$ is called a {\it coboundary} for
$\Cal G$ and the set of coboundaries is denoted $\text{\rm B}^1(\Cal G)$. It
is clearly a subgroup of $\text{\rm Z}^1(\Cal G)$. We denote the quotient
group $\text{\rm H}^1(\Cal G) \overset \text{\rm def} \to = \text{\rm Z}^1(\Cal
G)/\text{\rm B}^1(\Cal G)=\text{\rm Z}^1(\Cal G)/\sim_c$
and call it the $1${\it 'st
cohomology group} of $\Cal G$.

By the correspondence between countably generated full groupoids
and countable m.p. standard equivalence relations described in
Section 1.3, one can alternatively view the 1-cohomology groups
$\text{\rm Z}^1(\Cal G), \text{\rm B}^1(\Cal G), \text{\rm H}^1(\Cal G)$
as associated to the
equivalence relation $\Cal R=\Cal R_{\Cal G}$, in which case one
recovers the definition of $\text{\rm H}^1(\Cal R)$ from (page 308 of [FM]).

Let now $A \subset M$ be a II$_1$ factor with a Cartan subalgebra.
If $\theta \in {\text{\rm Aut}}_0(M;A)$ and $\phi_v = {\text{\rm
Ad}}(v)\in \Cal G_{A \subset M}$ for some $v \in \Cal G\Cal
N_M(A)$ then $w^\theta({\phi_v}) = \theta(v)v^*$ is a well defined
$1$-cocycle for $\Cal G$. Conversely, if $w\in \text{\rm H}^1(\Cal G)$ then
there exists a unique automorphism $\theta^w \in {\text{\rm
Aut}}_0(M;A)$ satisfying $\theta^w(a v) aw_{\phi_v}v$, $\forall
a\in A, v \in \Cal G\Cal N_M(A)$.

\proclaim{1.5.1. Proposition} $\theta \mapsto w^\theta$ is an
isomorphism of topological groups, from ${\text{\rm Aut}}_0(M;A)$
onto $\text{\rm Z}^1(\Cal G_{A \subset M})$, that takes ${\text{\rm
Int}}_0(M;A)=\{\text{\rm Ad}(u) \mid u \in \Cal U(A)\}$ onto
$\text{\rm B}^1(\Cal G_{A \subset M})$ and whose inverse is $w \mapsto
\theta^w$. Thus, $\theta \mapsto w^\theta$ implements an
isomorphism between the topological groups ${\text{\rm
Out}}_0(M;A)={\text{\rm Aut}}_0(M;A)/{\text{\rm Int}}_0(M;A)$ and
$\text{\rm H}^1(\Cal G_{A \subset M})$.
\endproclaim
\vskip .05in {\it Proof}. This is trivial by the definitions.
\hfill Q.E.D.

\vskip .05in By a well known lemma of Connes (see e.g. [C75]), if
$\theta \in {\text{\rm Aut}}_0(M;A)$ satisfies
$\theta_{|pMp}={\text{\rm Ad}}(u)_{|pMp}$ for some $p \in \Cal
P(A)$, $u \in \Cal U(A)$ then $\theta \in {\text{\rm
Int}}_0(M;A)$. Thus, $\theta \mapsto \theta_{|pMp}$ defines an
isomorphism from Out$_0(M;A)$ onto Out$_0(pMp;Ap)$. Applying this
to the Cartan subalgebra inclusion $L(\Cal G_0) \subset L(\Cal G)$
for $\Cal G$ an abstract ergodic full groupoid acting on the
non-atomic probability space, from 1.5.1 we get: $\text{\rm
H}^1(\Cal G)$ is naturally isomorphic to $\text{\rm H}^1(\Cal
G^t), \forall t > 0$. In particular, since 1.5.1 also implies
$\text{\rm H}^1(\sigma)=\text{\rm H}^1(\Cal G_\sigma)$, it follows
that $\text{\rm H}^1(\sigma)$ is invariant to stable orbit
equivalence. We have thus shown:

\proclaim{1.5.2. Corollary} $1^\circ$. $\text{\rm H}^1(\Cal G^t)$ is
naturally isomorphic to $\text{\rm H}^1(\Cal G), \forall t > 0$.

$2^\circ$. If $\sigma$ is a free ergodic measure preserving action
then $\text{\rm H}^1(\sigma) = \text{\rm H}^1(\Cal G_\sigma)$ and
$\text{\rm H}^1(\sigma)$ is
invariant to stable orbit equivalence. Also, $\text{\rm
Z}^1(\sigma)=\text{\rm Z}^1(\Cal G_\sigma)$ and $\text{\rm
Z}^1(\sigma)$ is invariant to orbit equivalence.
\endproclaim

Note that the equality $\text{\rm H}^1(\sigma) = \text{\rm
H}^1(\Cal G_\sigma)$ (and thus the  invariance of $\text{\rm
H}^1(\sigma)$ to orbit equivalence) was already shown in ([FM77]).

\vskip .05in \noindent {\it 1.6. The closure of $\text{\rm
B}^1(\Cal G)$ in $\text{\rm Z}^1(\Cal G)$}. Given any ergodic full
groupoid $\Cal G$, the groups $\text{\rm B}^1(\Cal
G)\simeq{\text{\rm Int}}_0(M;A)$ are naturally isomorphic to $\Cal
U(A)/\Bbb T$, where $A = L(\Cal G_0), M=L(\Cal G)$. But this
isomorphism doesn't always carry the topology that $\text{\rm
B}^1(\Cal G)$ (resp. ${\text{\rm Int}}_0(M;A)$) inherit from
$\text{\rm Z}^1(\Cal G)$ (resp. Aut$_0(M;A)$) onto the quotient of
the $\| \cdot\|_2$-topology on $\Cal U(A)/\Bbb T$. It was shown by
K. Schmidt in ([S80], [S81]) that the two topologies on $\text{\rm
B}^1(\sigma)$ coincide iff the action $\sigma$ is strongly
ergodic. We recall his result in the statement below, relating it
to a result of A. Connes, showing that the group of inner
automorphisms of a II$_1$ factor is closed iff the factor has no
non-trivial central sequences ([C75]).

\proclaim{1.6.1. Proposition} Let $A \subset M$ be a II$_1$ factor
with a Cartan subalgebra. The following conditions are equivalent:

$(a)$. $\text{\rm H}^1(\Cal G_{A\subset M})$ is a Polish group (equivalently
$\text{\rm H}^1(\Cal G_{A \subset M})$
is separate), i.e. $\text{\rm B}^1(\Cal G_{A
\subset M})$ is closed in $\text{\rm Z}^1(\Cal G_{A \subset M})$.

$(b)$. ${\text{\rm Int}}_0(M;A)$ is closed in ${\text{\rm
Aut}}_0(M;A)$.

$(c).$ The action of $\Cal G_{A \subset M}$ on $A$ is strongly
ergodic, i.e. it has no non-trivial asymptotically invariant
sequences.

$(d)$. $M'\cap A^\omega = \Bbb C$, where $\omega$ is a free
ultrafilter on $\Bbb N$.

Moreover, if $M = A \rtimes_\sigma G$ for some free action $\sigma$ of a
group $G$ on $(A, \tau)$, then the above conditions are equivalent
to $\sigma$ being strongly ergodic.
\endproclaim
\noindent {\it Proof}. $(a) \Leftrightarrow (b)$ follows from
1.5.1 and  $(c) \Leftrightarrow (d)$ is well known (and trivial).
Then notice that  $(b) \Leftrightarrow (d)$ is a relative version
of Connes' result in ([C75]), showing that ``${\text{\rm Int}}(N)$
is closed in Aut$(N)$ iff $N$ has no non-trivial central
sequences'' for II$_1$ factors $N$. Thus, a proof of $(b)
\Leftrightarrow (d)$ is obtained by following the argument in
([C75]), but replacing everywhere Int$(N)$ by Int$_0(M;A)$,
Aut$(N)$ by Aut$_0(M;A)$ and ``non-trivial central sequences of
$N$'' by ``non-trivial central sequences of $M$ that are contained
in $A$''.

To prove the last part, note that $\sigma$ strongly ergodic iff
$\{u_g\}_g'\cap A^\omega = \Bbb C$, where $\{u_g\}_g \subset M$
denote the canonical unitaries implementing the action $\sigma$ of
$G$ on $A$. But $\{u_g\}_g' \cap A^\omega = (A \cup \{u_g\}_g)'
\cap A^\omega = M'\cap A^\omega$, hence strong ergodicity of
$\sigma$ is equivalent to $(d)$. \hfill Q.E.D.

\vskip .05in It was shown in ([S80], [S81]) that arbitrary ergodic
m.p. actions $\sigma$ of infinite property (T) groups $G$ are
always strongly ergodic, and that in fact $\text{\rm B}^1(\sigma)$
is always open (and thus also closed) in $\text{\rm Z}^1(\sigma)$.
The interpretation of the inclusion $\text{\rm B}^1(\sigma)
\subset \text{\rm Z}^1(\sigma)$ as $\text{\rm Int}_0(M;A) \subset
\text{\rm Aut}_0(M;A)$ makes this result into a relative version
of the rigidity result in ([C80]), showing that for property (T)
factors Int$(M)$ is open in Aut$(M)$. We notice here the following
generalization:

\proclaim{1.6.2. Proposition} Assume $G$ has an infinite subgroup
$H \subset G$ such that the pair $(G,H)$ has the relative property
$(\text{\rm T})$. If $\sigma$ is a free m.p. action of $G$ on the
probability space such that $\sigma_{|H}$ is ergodic, then
$\sigma$ is strongly ergodic, equivalently
$\text{\rm B}^1(\sigma)$ is closed
in $\text{\rm Z}^1(\sigma)$. Moreover,
the subgroup $\text{\rm Z}^1_H(\sigma)\overset
\text{\rm def} \to = \{w\in \text{\rm Z}^1(\sigma) \mid w_{|H} \sim_c
{\bold 1}_H\}$ is open and closed in $\text{\rm Z}^1(\sigma)$.
\endproclaim
\noindent {\it Proof}. Since $(G,H)$ has the relative property
(T), by ([Jol02]) there exist a finite subset $F \subset G$ and
$\delta > 0$ such that if $\pi: G \rightarrow \Cal U(\Cal H)$,
$\xi\in \Cal H$, $\|\xi\|_2=1$ satisfy $\|\pi_g(\xi)-\xi\|_2 \leq
\delta, \forall g \in F$ then $\|\pi_h(\xi)-\xi\|_2 \leq 1/2,
\forall h\in H$ and $\pi_{|H}$ has a non-trivial fixed vector.

If $\sigma$ is not strongly ergodic then there exists $p \in \Cal
P(A)$ such that $\tau(p)=1/2$ and $\|\sigma_g(p)-p\|_2 \leq
\delta/2, \forall g\in F$. But then $u=1-2p$ satisfies $\tau(u)=0$
and $\|\sigma_g(u)-u\|_2 \leq \delta, \forall g\in F$. Taking
$\pi$ to be the $G$-representation induced by $\sigma$ on $L^2(A,
\tau) \ominus \Bbb C1$, it follows that $L^2(A, \tau) \ominus \Bbb
C1$ contains a non-trivial vector fixed by $\sigma_{|H}$. But this
contradicts the ergodicity of $\sigma_{|H}$.

Let now $M=A \rtimes_\sigma G$ and $\theta=\theta^w \in {\text{\rm
Aut}}_0(M;A)$ be the automorphism associated to some $w\in
\text{\rm Z}^1(\sigma)$ satisfying $\|\theta(u_g)-u_g\|_2=\|w_g-1\|_2 \leq
\delta, \forall g\in F$. Then the unitary representation $\pi: G
\rightarrow \Cal U(L^2(M, \tau))$ defined by $\pi_g(\xi)=u_g \xi
\theta(u_g^*)$ satisfies
$\|\pi_g(\hat{1})-\hat{1}\|_2=\|w_g-1\|_2\leq \delta$, $\forall
g\in F$. Thus, $\|w_h-1\|_2 = \|\pi_h(\hat{1})-\hat{1}\|_2 \leq 1/2$
implying
$$
\|\theta(vu_h)-vu_h\|_2=\|\theta(u_h)-u_h\|_2 \leq 1/2, \forall
h\in H, v\in \Cal U(A).
$$
It follows that if $b$ denotes the element of minimal norm
$\|\cdot \|_2$ in $\overline{\text{\rm co}}^w \{u_h^*v^*
\theta(vu_h) \mid h \in H, v\in \Cal U(A) \}$ then $\|b-1\|_2 \leq
1/2$ and $vu_h b = b\theta(vu_h) = b w_h vu_h$, $\forall h \in H,
v \in \Cal U(A)$. But this implies $b \neq 0$ and $xb =b \theta(x),
\forall x\in N=A \rtimes_{\sigma_{|H}} H$. In particular
$[b,A]=0$ so $b \in A\subset N $. Since $N$ is a factor (because
$\sigma_{|H}$ is ergodic), this implies $b$ is a scalar multiple
of a unitary element $u$ in $A$ satisfying $w_h=u^*\sigma_h(u),
\forall h\in H$. Thus $w\in \text{\rm Z}^1_H(\sigma)$, showing that
$\text{\rm Z}^1_H(\sigma)$ is open (thus closed too). \hfill Q.E.D.

\proclaim{1.6.3. Corollary} If $G$ has an infinite subgroup $H
\subset G$ such that the pair $(G,H)$ has the relative property
$(\text{\rm T})$ then $\text{\rm Char}(G)/\{\gamma \in \text{\rm
Char}(G) \mid \gamma_{|H}={\bold 1}_H\}$ is a finite group.
\endproclaim
\noindent {\it Proof}. Take $\sigma$ to be the (classic)
$G$-Bernoulli shift and view $\text{\rm Char}(G)$ as a subgroup of
$Z^1(\sigma)$, as in 1.1. By 1.6.2 it follows that $\gamma \in
\text{\rm Char}(G)$ satisfies $\gamma_{|H}\sim_c {\bold 1}_H$ iff
$\gamma_{|H}={\bold 1}_H$ (because $\sigma_{|H}$ is mixing; see
also 2.4.1$^\circ$). Thus, Z$^1_H(\sigma) \cap \text{\rm Char}(G)=
\{\gamma \in \text{\rm Char}(G) \mid \gamma_{|H}={\bold 1}_H\}$
and the statement follows from 1.6.2$^\circ$. \hfill Q.E.D.

\vskip .05in Other actions shown to be strongly ergodic in ([S80],
[S81]) are Bernoulli shifts of non-amenable groups and the action
of $SL(2, \Bbb Z)$ on $(\Bbb T^2, \lambda)$. There is in fact a
common explanation for these examples: In both cases the
representation implemented by $(\sigma, G)$ on $L^2(X, \mu)
\ominus \Bbb C1$ can be realized as $\oplus_i \ell^2(G/\Gamma_i)$
with $\Gamma_i \subset G$ amenable subgroups. (If $\sigma$ is a
$G$-Bernoulli shift, then $\Gamma_i$ can even be taken finite, see
also [J83b]. If in turn $G=SL(2, \Bbb Z)$ then the action $\sigma$
of $G$ on $L^2(\Bbb T^2, \lambda)\ominus \Bbb C1 \simeq
\ell^2(\Bbb Z^2 \setminus \{(0,0)\})$ corresponds to the action of
$G$ on $\Bbb Z^2\setminus \{(0,0)\}$ and $\Gamma_i$ are
stabilizers of elements $h_i \in \Bbb Z^2 \setminus \{(0,0)\}$,
thus amenable.) Hence, if $(\sigma, G)$ would not be strongly
ergodic then the trivial representation of $G$ would be weakly
contained $\oplus_i \ell^2(G/H_i)$ and the following general
observation applies:

\proclaim{1.6.4. Lemma}  Let $G$ be a non-amenable group and
$\{H_i\}_i$ the family of amenable subgroups of $G$. Then the
trivial representation of $G$ is not weakly contained in $\oplus_i
\ell^2(G/H_i)$ (thus not weakly contained in $\oplus_i
\ell^2(G/H_i)\overline{\otimes} \ell^2(\Bbb N)$ either).
\endproclaim
\noindent {\it Proof}. This follows immediately from the
continuity of induction of representations. Indeed, every
$\ell^2(G/H_i)$ is equivalent to the induced from $H_i$ to $G$ of
the trivial representation $1_{H_i}$ of $H_i$, Ind$_{H_i}^G
1_{H_i}.$ Since $H_i$ is amenable, $1_{H_i}$ follows weakly
contained in the left regular representation $\lambda_{H_i}$ of
$H_i.$ Thus, Ind$_{H_i}^G 1_{H_i}$ is weakly contained in
Ind$_{H_i}^G (\lambda_{H_i})$, which in turn is just the left
regular representation $\lambda_G$ of $G.$ Altogether, this shows
that if $1_G$ is weakly contained in $\oplus_i \ell^2(G/H_i)$ then
it is weakly contained in a multiple of $\lambda_G$. Since the
latter is weakly equivalent to $\lambda_G$, $1_G$ follows weakly
contained in $\lambda_G$, implying that $G$ is amenable, a
contradiction. \hfill Q.E.D.

\vskip .05in If one defines the property $(\tau)$ for a group $G$
with respect to a family $\Cal L$ of subgroups by requiring that
the trivial representation of $G$ is not an accumulation point of
$\oplus_{H \in \Cal L} \ell^2(G/H)$, like in ([L94]), then the
above Lemma can be restated: {\it If $G$ is non-amenable then it
has the property $(\tau)$ with respect to the family $\Cal L$ of
its amenable subgroups}.

Let now $G$ act by automorphisms on a discrete group $H$ and
denote by $\sigma$ the action it implements on the finite group
von Neumann algebra $(L(H), \tau)$, then note that the ensuing
representation of $G$ on $L^2(L(H)\ominus \Bbb C1)=\ell^2(H
\setminus \{e\})$ is equal to $\oplus_h \ell^2(G/\Gamma_h)$, where
$\Gamma_h \subset G$ denotes the stabilizer of $h\in H \setminus
\{e\}$, $\Gamma_h=\{\gamma \in G \mid \gamma(h)=h\}$. Lemma 1.6.4
thus shows:

\proclaim{1.6.5. Corollary} Assume $G$ is non-amenable and the
stabilizer of each $h \in H\setminus \{e\}$ is amenable. For any
non-amenable $\Gamma \subset G$ the action $\sigma_{|\Gamma}$ of
$\Gamma$ on $(L(H), \tau)$ is strongly ergodic. In particular, if
$\Gamma \subset SL(2, \Bbb Z)$ is non-amenable then the
restriction to $\Gamma$ of the canonical action of $SL(2, \Bbb Z)$
on $(\Bbb T^2, \mu)$ is strongly ergodic. \endproclaim

Finally, note that if one takes $H=G$ and let $G$ act on itself by
conjugation, then Lemma 1.6.4 implies that if $G$ is non-amenable
and the commutant in $G$
of any $h \in G \setminus \{e\}$ is amenable then $G$ is not inner
amenable either.

\heading 2. 1-cohomology for quotients of $G$-Bernoulli shifts
\endheading

In this section we consider groups $G$ satisfying some ``mild''
rigidity property and construct many examples of actions
$(\sigma', G)$ for which we can explicitly calculate the
1-cohomology. The actions $\sigma'$ are quotients of the
$G$-Bernoulli shift $\sigma$ (or of other ``malleable'' actions of
$G$), obtained by restricting $\sigma$ to subalgebras that are
fixed points of groups of automorphisms in the commutant of
$\sigma$. The construction is inspired from ([P01]), where a
similar idea is used to produce actions on the hyperfinite II$_1$
factor that have prescribed fundamental group and prescribed
1-cohomology.

This calculation of $\text{\rm H}^1(\sigma',G)$ works whenever the
1-cohomology group of the ``initial'' action $(\sigma,G)$ is equal
to the character group of $G$. For $G$ weakly rigid and $\sigma$
Bernoulli shifts, $\text{\rm H}^1(\sigma,G)$ was shown equal to
Char$(G)$ in ([PSa03]), by adapting to the commutative case the
proof of the similar result for non-commutative Bernoulli shifts
in ([P01]). We begin by extending the result in ([PSa03]) to a
more general context in which the argument in ([P01], [PSa03])
still works.

\vskip .05in \noindent {\it 2.1. Definition} ([P01], [P03]). An
integral preserving action $\sigma : G \rightarrow {\text{\rm
Aut}}(A, \tau)$ of $G$ on $A\simeq L^\infty(X, \mu)$ is {\it
w-malleable} if there exist a decreasing sequence of abelian von
Neumann algebras $\{(A_n, \tau)\}_n$ containing $A$ and actions
$\sigma_n : G \rightarrow {\text{\rm Aut}}(A_n, \tau)$, such that
$\cap A_n = A$, ${\sigma_n}_{|A_{n+1}}=\sigma_{n+1}$,
${\sigma_n}_{|A}=\sigma$, $\forall n$ and such that for each $n$
the flip automorphism $\alpha_1$ on $A_n \overline{\otimes} A_n$,
definded by $\alpha_1(x \otimes y) = y \otimes x$, $x,y\in A_n$,
is in the connected component of the identity in the Polish group
$\tilde{\sigma_n}(G)' \cap {\text{\rm Aut}}(A_n \overline{\otimes}
A_n, \tau \times \tau)$, where $\tilde{\sigma_n}$ is the
automorphism on $A_n \overline{\otimes} A_n$ given by
$\tilde{\sigma}_n(g)=\sigma_n (g) \otimes \sigma_n(g), g\in G$. If
$H \subset G$ is a subgroup, then the action $\sigma$ is {\it
w-malleable} {\it w-mixing}$/H$ (resp. {\it w-malleable mixing})
if the extensions $\sigma_n$ can be chosen so that
${\sigma_n}_{|H}$ are weakly mixing (resp. so that $\sigma_n$ are
mixing), $\forall n$.

\vskip .05in \noindent {\it 2.1' Example}. Let $(g,s) \mapsto gs$
be an action of the group $G$ on a set $S$, $(Y_0, \nu_0)$ be a
non-trivial standard probability space. Let $(X, \mu)=\Pi_s (X_0,
\nu_0)_s$ and denote by $\sigma$ the {\it Bernoulli shift action}
of $G$ on $L^\infty(X, \mu)$ implemented by
$\sigma_g((x_s)_s)=(x'_s)_s$, where $x'_s = x_{gs}$. If $(g,s)
\mapsto gs$ satisfies: \vskip .05in \noindent $(2.1.1)$. $\forall
g \neq e, \exists s\in S$ with $gs\neq s$, \vskip .05in \noindent
then $\sigma$ is free. If $H \subset G$ is a subgroup such that:
\vskip .05in \noindent $(2.1.2)$. $\forall S_0 \subset S$ finite
$\exists F_\infty \subset H$ infinite set with $hS_0 \cap S_0 =
\emptyset$, $\forall h \in F_\infty$, \vskip .05in \noindent then
$\sigma_{|H}$ is weakly mixing. Also, if: \vskip .05in \noindent
$(2.1.2')$. $\forall S_0 \subset S$ finite $\exists F_0 \subset G$
finite set with $hS_0 \cap S_0 = \emptyset$, $\forall h \in G
\setminus F_0$, \vskip .05in \noindent then $\sigma$ is mixing. If
$S=G$ and we let $G$ act on itself by left multiplication then
$\sigma$ is called a {\it classic Bernoulli shift} action of $G$.

\proclaim{2.2. Lemma} Let $G$ be a group with an infinite subgroup
$H \subset G$. A Bernoulli shift action $(\sigma, G)$ satisfying
$(2.1.1)$, $(2.1.2)$ is free, w-malleable w-mixing$/H$. Also, a
classic Bernoulli shift is free, w-malleable mixing (on $G$).
\endproclaim
\vskip .05in \noindent {\it Proof}. The proof of (3.2 in [PoSa])
shows that if $(Y_0, \nu_0) \simeq (\Bbb T, \lambda)$  then
$\sigma$ is malleable. For general $(Y_0, \nu_0)$, the proof of
(3.6 in [PoS03]) shows that the action $\sigma$ can be
``approximated from above'' by Bernoulli shifts with base space
$\simeq (\Bbb T, \lambda)$, all satisfying $(2.1.2)$, thus being
w-malleable w-mixing$/H$. If $\sigma$ is a classic Bernoulli shift
then $(2.1.2')$ is satisfied, so $\sigma$ is mixing. \hfill Q.E.D.

\vskip .05in \noindent {\it 2.3. Definition}. Let $G$ be a group.
An infinite subgroup $H\subset G$ is {\it wq-normal} in $G$ if
there exists a countable ordinal $\imath$ and a well ordered
family of intermediate subgroups $H = H_0 \subset H_1 \subset ...
\subset H_{\jmath} \subset ... \subset H_{\imath} = G$ such that
for each $\jmath < \imath$, $H_{\jmath+1}$ is the group generated
by the elements $g \in G$ with $gH_{\jmath}g^{-1} \cap H_{\jmath}$
infinite and such that if $\jmath \leq \imath$ has no
``predecessor'' then $H_{\jmath} = \cup_{n < \jmath} H_n$. Note
that this condition is equivalent to the following: \vskip .05in
\noindent $(2.3')$. There exists no intermediate subgroup $H
\subset K \varsubsetneq G$ such that $gKg^{-1} \cap K$ is finite
$\forall g\in G\setminus K$. Equivalently, for all $H\subset K
\varsubsetneq G$ there exists $g \in G\setminus K$, $gKg^{-1} \cap
K$ is infinite.

\vskip .05in

Indeed, if $H \subset G$ satisfies $(2.3')$
then it clearly satisfies $2.3$, by (countable transfinite) induction.
Conversely, assume $H \subset G$ satisfies
$2.3$ and let $K \varsubsetneq G$ be a subgroup containing
$H$ such that $gKg^{-1} \cap K$
finite $\forall g\in G\setminus K$. We show that
this implies $K=G$, giving a contradiction. It is sufficient to show
that $H_{\jmath}\subset K$ implies $H_{\jmath +1} \subset K$. If
there exists $g\in H_{\jmath+1} \setminus K=  H_{\jmath+1} \setminus K
\cap  H_{\jmath+1}$ then we would have $gKg^{-1}\cap K$ finite so in
particular $gH_{\jmath}g^{-1} \cap H_{\jmath}$ finite. But this implies that
all $g \in H_{\jmath+1}$ for which $gH_{\jmath}g^{-1}
\cap H_{\jmath}$ is infinite lie
in $K \cap H_{\jmath +1}$, thus $K \supset H_{\jmath+1}$ by the way
$H_{\jmath+1}$ was defined.

By condition 2.3 we see that if $H\subset G$ is wq-normal and $G$
is embedded as a normal subgroup in some larger group
$\overline{G}$ (or even merely as a {\it quasi-normal} subgroup
$G\subset \overline{G}$, i.e. so that $gGg^{-1} \cap G$ has finite
index in $G$, $\forall g\in \overline{G}$) then $H \subset
\overline{G}$ is wq-normal. Condition $(2.3')$ shows that an
inclusion of groups of the form $H \subset G=H* H'$, with $H$
infinite and $H'$ non-trivial is not wq-normal.

\proclaim{2.4. Lemma} Let $G$ be an infinite group and $\sigma$ a
free ergodic measure preserving action of $G$ on the probability
space.

$1^\circ$. Assume $\sigma$ is weakly mixing and for each $\gamma
\in {\text{\rm Char}}(G)$ denote $w^\gamma$ the $1$-cocycle
$w^\gamma_g=\gamma(g)1, g\in G$. Then the group morphism $\gamma
\mapsto w^\gamma$ is $1$ to $1$ and continuous from $\text{\rm
Char}(G)$ into $\text{\rm H}^1(\sigma,G)$.

$2^\circ$. Assume $H \subset G$ is an infinite subgroup of $G$
such that either $H$ is normal in $G$ and $\sigma_{|H}$ is weakly
mixing or $H$ is weakly normal in $G$ and $\sigma$ is mixing. If
$w \in \text{\rm Z}^1(\sigma, G)$ is so that $w_{|H}\in {\text{\rm Char}}(H)$
then $w \in {\text{\rm Char}}(G)$.
\endproclaim
\vskip .05in \noindent
{\it Proof}. $1^\circ$. If $w_1(g) = u^*
w_2(g)\sigma_g(u), \forall g\in G$ then $\sigma_g(u) \in \Bbb C u,
\forall g\in G$ and since $\sigma$ is weakly mixing, this implies
$u \in \Bbb C1$ so $w_1 = w_2$.

$2^\circ$. In both cases, it is clearly sufficient to prove that
if $g_0 \in G$ is so that $H'=g_0^{-1}Hg_0 \cap H$ is infinite and
$\sigma$ is weakly mixing on $H'$ with $w_{|H}=\gamma \in
{\text{\rm Char}}(H)$ then $w_{g_0} \in \Bbb C1$. To see this, let
$M = A \rtimes_\sigma G$ with $\{u_g\}_g\subset M$ denoting the
canonical unitaries implementing the action $\sigma$ of $G$ on
$A$. Also, denote $u'_{g}=w_{g}u_{g}, g\in G$. If we let $h'\in
H'$ and denote $h=g_0h'g_0^{-1}$, then $h\in H$ and we have
$$
\gamma_{h'} w_{g_0} \sigma_{h}(w_{g_0}^*)u_{h}=\gamma_{h'} w_{g_0}
u_{h} w_{g_0}^*
$$
$$
=\gamma_{h'}w_{g_0} u_{g_0}u_{h'}u_{g_{0}^{-1}}w_{g_0}^*=
(w_{g_0}u_{g_0})(\gamma_{h'}u_{h'})(w_{g_0}u_{g_0})^*
$$
$$
=u'_{g_0}u'_{h'}u'_{{g_0}^{-1}}=u'_{g_0h'g_0^{-1}}=
w_{g_0h'g_0^{-1}}u_{g_0h'g_0^{-1}}=\gamma_{h}
u_{h}.
$$
Thus, $\sigma_{h}(w_{g_0}) \in \Bbb C w_{g_0}, \forall h \in
g_0H'g_0^{-1}$ and since $\sigma_{|g_0H'g_0^{-1}}$ is weakly
mixing (because $\sigma_{|H'}$ is weakly mixing) this implies
$w_{g_0} \in \Bbb C1$.

\hfill Q.E.D.

\proclaim{2.5. Corollary} Let $H \subset G$, $\sigma$ be as in
$2.4.2^\circ$. If the restriction to $H$ of any $w \in \text{\rm Z}^1(\sigma,
G)$ is cohomologous to a character of $H$ then $\text{\rm H}^1(\sigma,
G)=\text{\rm Char}(G)$.
\endproclaim

For the next statement, recall from ([Ma82]) that if $H \subset G$
is an inclusion of groups then the {\it pair} $(G,H)$ {\it has the
relative property} (T) if all unitary representations of $G$ that
weakly contain the trivial representation of $G$ must contain the
trivial representation of $H$ (when restricted to $H$).

\proclaim{2.6. Theorem} Let $G$ be a countable discrete group with
an infinite subgroup $H \subset G$ such that $(G,H)$ has the
relative property $(\text{\rm T})$. Let $\sigma$ be a free ergodic
m.p. action of $G$ on the probability space. Assume that either
$\sigma$ is w-malleable mixing or that it is w-malleable
w-mixing$/H$. Then any $1$-cocycle $w$ for $(\sigma,G)$ is
cohomologous to a $1$-cocycle which is scalar valued when
restricted to $H$, i.e. $\exists \gamma \in \text{\rm Char}(H)$
such that $w_{|H} \sim_c \gamma {\bold 1}_H$.  Thus, if in
addition $H$ is wq-normal in $G$ then $\text{\rm
H}^1(\sigma,G)={\text{\rm Char}}(G)$.
\endproclaim
\vskip .05in \noindent {\it Proof}. Let $w\in \text{\rm
Z}^1(\sigma,G)$. The proof that if $\sigma$ is a classic Bernoulli
shift (thus w-malleable mixing by Lemma 2.2) then $\text{\rm
H}^1(\sigma,G)= \text{\rm Char}(G)$ in ([PoS03]) only uses the
w-malleability of $\sigma$ to derive that $w_{|H}$ is cohomologous
to a character of $H$. But then Lemma 2.4 shows that $w$ is
cohomologous to a character of $G$, so $\text{\rm
H}^1(\sigma,G)={\text{\rm Char}}(G)$ by 2.5. \hfill Q.E.D.

\proclaim{2.7. Lemma} Let $G, \Gamma$ be discrete groups with $G$
infinite. Let $\sigma$ be a free, weakly mixing m.p. action of $G$
on the probability space and $\beta$ a free measure preserving
action of $\Gamma$ on the same probability space which commutes
with $\sigma$. If $A^\Gamma\overset \text{\rm def} \to =\{a\in
A\mid \beta_h(a)=a, \forall h\in \Gamma\}$ then
$\sigma_g(A^\Gamma)=A^\Gamma, \forall g\in G,$ so
$\sigma^\Gamma_g \overset \text{\rm def} \to =
{\sigma_g}_{|A^\Gamma}$ defines an integral preserving
action of $G$ on $A^\Gamma$.
\endproclaim
\vskip .05in \noindent {\it Proof}. Since $\beta_h(\sigma_g(a)) =
\sigma_g(\beta_h(a))=\sigma_g(a)$, $\forall h \in \Gamma, a\in
A^\Gamma$, it follows that $\sigma_g$ leaves $A^\Gamma$ invariant
$\forall g \in G$. \hfill Q.E.D.

\proclaim{2.8. Lemma} Under the same hypothesis and with the same
notations as in $2.7$, assume the action $\sigma^\Gamma$ of $G$ on
the fixed point algebra $A^{\Gamma}$ is free. For each $\gamma \in
\text{\rm Char}(\Gamma)$ denote $\Cal U_\gamma \overset \text{\rm
def} \to = \{v \in \Cal U(A)\mid \beta_h(v)=\gamma(h)v, \forall
h\in \Gamma\}$ and ${\text{\rm Char}}_\beta(\Gamma)\overset
\text{\rm def} \to = \{ \gamma \in \text{\rm Char}(\Gamma) \mid
\Cal U_\gamma \neq \emptyset \}$. Then we have:

$1^\circ$. $\Cal U_\gamma \Cal U_{\gamma'} = \Cal U_{\gamma
\gamma'}, \forall \gamma, \gamma' \in \text{\rm Char}(G)$, and $\text{\rm
Char}_\beta(\Gamma)$ is a countable group.

$2^\circ$. If $\gamma_0 \in \text{\rm Char}(G)$, $\gamma \in
\text{\rm Char}_\beta(\Gamma)$ and $v \in \Cal U_\gamma$, then
$w^{\gamma_0, \gamma}(g) \overset \text{\rm def} \to =
\sigma_g(v)v^* \gamma_0(g) \in A^\Gamma, \forall g\in G$, and
$w^{\gamma_0, \gamma}$ defines a $1$-cocycle for $(\sigma^\Gamma,
G)$ whose class in $\text{\rm H}^1(\sigma^\Gamma, G)$ doesn't depend on the
choice of $v \in \Cal U_\gamma$.
\endproclaim
\vskip .05in \noindent {\it Proof}. $1^\circ$. If $v\in \Cal
U_\gamma, v'\in \Cal U_{\gamma'}$ then
$$
\beta_h(v v')=\beta_h(v)\beta_h(v')
=\gamma(h)\gamma'(h) vv',
$$
so $vv'\in \Cal U_{\gamma \gamma'}$. This also implies $\text{\rm
Char}_\beta(\Gamma)$ is a group. Noticing that $\{\Cal
U_\gamma\}_\gamma$ are mutually orthogonal in
$L^2(A,\tau)=L^2(X,\mu)$, by the separability of $L^2(X,\mu)$,
$\text{\rm Char}_\beta(\Gamma)$ follows countable.

$2^\circ$. Since $\sigma, \beta$ commute, $\sigma_g(\Cal U_\gamma)
= \Cal U_\gamma$, $\forall g\in G, \gamma \in \text{\rm
Char}_\beta(\Gamma)$. In particular, $\sigma_g(v)v^*\in \Cal U_1=
\Cal U(A^\Gamma)$, $\forall g\in G$ showing that the function
$w^{\gamma_0, \gamma}$ takes values in $\Cal U(A^\Gamma)$. Since
$w^{\gamma_0, \gamma}$ is clearly a $1$-cocylce for $\sigma$ (in
fact $w^{\gamma_0, \gamma} \sim_c \gamma_01$ as elements in
$\text{\rm Z}^1(\sigma,G)$), it follows that $w^{\gamma_0, \gamma}\in
\text{\rm Z}^1(\sigma^\Gamma,G)$.

If $v'$ is another element in $\Cal U_\gamma$ then $u = v'v^*\in
\Cal U(A^\Gamma)$ and the associated 1-cocycles $w^{\gamma_0,
\gamma}$ constructed out of $v, v'$ follow cohomologous via $u$ in
$\text{\rm Z}^1(\sigma^\Gamma, G)$. \hfill Q.E.D.

\proclaim{2.9. Theorem} With the same assumptions and notations as
in $2.8$, if $\text{\rm Char}_\beta(\Gamma)$ is given the discrete
topology then $\Delta: \text{\rm Char}(G) \times \text{\rm
Char}_\beta(\Gamma)\rightarrow \text{\rm H}^1(\sigma^\Gamma, G)$ defined by
letting $\Delta(\gamma_0, \gamma)$ be the class of $w^{\gamma_0,
\gamma}$ in $\text{\rm H}^1(\sigma^\Gamma,G)$ is a $1$ to $1$ continuous
group morphism. If in addition $\text{\rm H}^1(\sigma, G)=\text{\rm Char}(G)$
then $\Delta$ is an isomorphism of topological groups.
\endproclaim
\vskip .05in \noindent {\it Proof}. The map $\Delta$ is clearly a
group morphism and continuous. To see that it is 1 to 1 let
$\gamma_0 \in \text{\rm Char}(G)$, $\gamma \in \text{\rm
Char}_\beta(\Gamma)$ and $v\in \Cal U_\gamma$ and represent the
element $\Delta(\gamma_0, \gamma)\in \text{\rm H}^1(\sigma^\Gamma,G)$ by the
1-cocycle $w^{\gamma_0, \gamma}_g = \sigma_g(v)v^* \gamma_0(g),
g\in G$. If $w^{\gamma_0, \gamma} \sim_c {\bold 1}$ then there
exists $u \in \Cal U(A^\Gamma)$ such that
$\sigma_g(u)u^*=\sigma_g(v)v^* \gamma_0(g), \forall g\in G$. Thus,
if we denote $u_0 = uv^* \in \Cal U(A)$ then
$\sigma_g(u_0)u_0^*=\gamma_0(g)1, \forall g$.
It follows that $\sigma_g(\Bbb Cu_0)=\Bbb Cu_0$,
$\forall g\in G$, and since $\sigma$ is weakly mixing
this implies $u_0 \in \Bbb C1$ and $\gamma_0=1$. Thus, $v\in \Bbb Cu \subset
\Cal U(A^{\Gamma})=\Cal U_1$, showing that $\gamma =1$ as well.

If we assume $\text{\rm H}^1(\sigma, G)=\text{\rm Char}(G)$ and
take $w\in \text{\rm Z}^1(\sigma^\Gamma,G)$ then we can view $w$
as a $1$-cocycle for $\sigma$. But then $w\sim_c \gamma_01$, for
some $\gamma_0 \in {\text{\rm Char}}(G)$. Since $\sigma$ is
ergodic, there exists a unique $v\in \Cal U(A)$ (up to
multiplication by a scalar) such that $w_g = \sigma_g(v)v^*
\gamma_0(g)$, $\forall g\in G$. Since $w$ is $A^\Gamma$-valued,
$\sigma_g(v)v^* \in \Cal U(A^\Gamma), \forall g$.  Thus
$\sigma_g(v)v^*=\beta_h(\sigma_g(v)v^*)=\sigma_g(\beta_h(v))
\beta_h(v)^*, \forall g$. By the uniqueness of $v$ this implies
$\beta_h(v)=\gamma(h)v$, for some scalar $\gamma(h)$. The map
$\Gamma \ni h \mapsto \gamma(h)$ is easily seen to be a character,
so $w=w^{\gamma_0, \gamma}$ showing that $(\gamma_0, \gamma)
\mapsto w^{\gamma_0, \gamma}$ is onto.

Since $\text{\rm H}^1(\sigma,G)=\text{\rm Char}(G)$ is compact, by
1.1 and 1.6.1 $\sigma$ is strongly ergodic so $\sigma^\Gamma$ is
also strongly ergodic. Thus $\text{\rm H}^1(\sigma^\Gamma,G)$ is
Polish, with $\Delta(\text{\rm Char}(G))$ a closed subgroup,
implying that $\Delta(\text{\rm Char}_\beta(\Gamma))\simeq
\text{\rm H}^1(\sigma^\Gamma)/\Delta(\text{\rm Char}(G))$ is
Polish. Since it is also countable, it is discrete. Thus, $\Delta$
is an isomorphism of topological groups. \hfill Q.E.D.

\vskip .05in
Note that in the above proof, from the hypothesis
H$^1(\sigma,G)= \text{\rm Char}(G)$ we only used the following
fact:

\vskip .05in \noindent $(2.9')$. There exists a continuous group
morphism $\text{\rm H}^1(\sigma,G) \ni \hat{w}_0 \mapsto
w_0 \in \text{\rm Z}^1(\sigma,G)$ retract of the quotient map $\text{\rm
Z}^1(\sigma,G) \rightarrow \text{\rm H}^1(\sigma,G)$ such that
each $w_0$ is $A^\Gamma$-valued (so that it can be viewed as an
element $w_0\in \text{\rm Z}^1(\sigma^\Gamma,G)$). \vskip .05in

Thus, if the condition $(2.9')$ is satisfied then the above proof
of $2.9$ shows that $\text{\rm H}^1(\sigma^\Gamma,G) \simeq
\text{\rm H}^1(\sigma,G) \times \text{\rm Char}_\beta (\Gamma)$.

\proclaim{2.10. Lemma} Let $G$ be an infinite group and $\sigma$
be the Bernoulli shift action of $G$ on $(X, \mu)=\Pi_g(\Bbb T,
\lambda)_g$. With the notations of $2.8, 2.9$, for any countable
abelian group $\Lambda$ there exists a countable abelian group
$\Gamma$ and a free action $\beta$ of $\Gamma$ on $(X, \mu)$ such
that $\text{\rm Char}_\beta(\Gamma)=\Lambda$, $[\sigma, \beta]=0$
and $\sigma_{|A^\Gamma}$ is a free action of $G$. Moreover, if
$\Lambda$ is finite then one can take $\Gamma=\Lambda$ and $\beta$
to be any action of $\Gamma=\Lambda$ on $(X,\mu)$ that commutes
with $\sigma$ and such that $\sigma \times \beta$ is a free action
of $G \times \Gamma$.
\endproclaim
\vskip .05in \noindent {\it Proof}. Let $\Gamma$ be a countable
dense subgroup in the (2'nd countable) compact group
$\hat{\Lambda}$ and $\mu_0$ be the Haar measure on
$\hat{\Lambda}$. Let $\beta_0$ denote the action of $\Gamma$ on
$L^\infty(\hat{\Lambda}, \mu_0)=L(\Lambda)$ given by
$\beta_0(h)(u_\gamma)=\gamma(h)u_\gamma, \forall h\in \Gamma$,
where $\{u_\gamma\}_{\gamma \in \Lambda} \subset L(\Lambda)$
denotes the canonical basis of unitaries in the group von Neumann
algebra $L(\Lambda)$ and $\gamma \in \Lambda$ is viewed as a
character on $\Gamma \subset \hat{\Lambda}$. Denote $A_0 =
L^\infty(\hat{\Lambda}, \mu_0) \overline{\otimes} L^\infty(\Bbb T,
\lambda)$ and $\tau_0$ the state on $A_0$ given by the product measure
$\mu_0 \times \lambda$. Let $\beta$ denote the product
action of $\Gamma$ on $\overline{\otimes}_{g\in G} (A_0,
\tau_0)_g$ given by $\beta(h) = \otimes_g(\beta_0(h) \otimes
id)_g$.

Since $(A_0, \tau_0) \simeq (L^\infty(\Bbb T, \lambda), \int \cdot
{\text{\rm d}}\lambda)$, we can view $\sigma$ as the Bernoulli
shift action of $G$ on $A=\overline{\otimes}_g (A_0, \tau_0)_g$.
By the construction of $\beta$ we have $[\sigma, \beta]=0$. Also,
the fixed point algebra $A^\Gamma$ contains a $\sigma$-invariant
subalgebra on which $\sigma$ acts as the (classic) Bernoulli shift. Thus,
the restriction $\sigma^\Gamma=\sigma_{|A^\Gamma}$ is a free,
mixing action of $G$. Finally, we see by construction that
$\text{\rm Char}_\beta(\Gamma)=\Lambda$.

The last part is trivial, once we notice that if the action
$\sigma \times \beta$ of $G \times \Gamma$ on $A$ is free then the
action $\sigma^\Gamma$ of $G$ on $A^\Gamma$ is free.
\hfill Q.E.D.

\vskip .05in From now on, it will be convenient to use the
following: \vskip .05in \noindent {\it 2.11. Notation}. We denote
by $w\Cal T$ the class of discrete countable groups $G$ which have
infinite, wq-normal subgroups $H\subset G$ such that the pair
$(G,H)$ has the relative property (T).

\vskip .05in

Note that all infinite property (T) groups are in the class $w\Cal
T$. Also, by 2.3 it follows that $w\Cal T$ is closed to inductive
limits and normal extensions (i.e. if $G \in w\Cal T$ and $G
\subset \overline{G}$ is a normal inclusion of groups then
$\overline{G}\in w\Cal T$). In particular, if $G\in w\Cal T$ and
$K$ is a group acting on $G$ by automorphisms then $G \rtimes K
\in w\Cal T$. For instance, if $G$ is infinite with property (T)
and $K$ is an arbitrary group then $G \times K \in w\Cal T$. Other
examples of groups in the class $w\Cal T$ are $\Bbb Z^2 \rtimes
SL(2,\Bbb Z)$ ([K67], [Ma82]), and more generally $\Bbb Z^2
\rtimes \Gamma$ for $\Gamma\subset SL(2,\Bbb Z)$ non-amenable (cf.
[B91]).

\proclaim{2.12. Corollary} Let $G \in w\Cal T$. Given any countable
discrete abelian group $\Lambda$ there exists a free ergodic m.p.
action $\sigma_\Lambda$ of $G$ on the standard non-atomic
probability space such that $\text{\rm H}^1(\sigma_\Lambda,G)=\text{\rm
Char}(G) \times \Lambda$. Moreover, if $\sigma$ denotes the
Bernoulli shift action of $G$ on $(X, \mu)=\Pi_g (\Bbb T, \mu)_g$
then all $\sigma_\Lambda$ can be taken to be quotients of
$(\sigma, (X, \mu))$ and such that the exact sequences of
$1$-cohomology groups $\text{\rm B}^1(\sigma_\Lambda) \hookrightarrow
\text{\rm Z}^1(\sigma_\Lambda) \rightarrow
\text{\rm H}^1(\sigma_\Lambda) \rightarrow 1$
are split. Thus, $\text{\rm Z}^1(\sigma_\Lambda)\simeq
\text{\rm H}^1(\sigma_\Lambda)
\times \text{\rm B}^1(\sigma_\Lambda)$.
\endproclaim
\vskip .05in \noindent {\it Proof}.  Since Bernoulli shifts with
base space $(\Bbb T, \mu)$ are malleable mixing, by ([PoS03]) or
Theorem 2.6 above we have $\text{\rm H}^1(\sigma,G)=\text{\rm
Char}(G)$ and the statement follows by Lemma 2.10 and Theorem 2.9.
The fact that $\sigma_\Lambda$ can be constructed so that its
exact sequence of 1-cohomology groups is split is clear from the
construction in the proof of 2.10, which shows that one can select
$u_\gamma \in \Cal U_\gamma$ such that $u_\gamma u_{\gamma'} =
u_{\gamma\gamma'}, \forall \gamma, \gamma' \in \Lambda =
{\text{\rm Char}}_\beta(\Gamma)$. \hfill Q.E.D.

\proclaim{2.13. Corollary} If $G\in w\Cal T$ then $G$ has a
continuous family of mutually non-stably orbit equivalent free
ergodic m.p. actions on the probability space, indexed by the
classes of virtual isomorphism of all
countable, discrete, abelian groups.
\endproclaim
\vskip .05in \noindent {\it Proof}. If we denote $K = \text{\rm
Char}(G)$ then $K$ is compact and open in $K \times \Lambda$.
Thus, for any isomorphism $\theta : K \times \Lambda_1 \simeq K
\times \Lambda_2$, $\theta(K)\cap K$ has finite index both in $K$
and in $\theta(K)$. Thus, $\Lambda_1, \Lambda_2$ must be virtually
isomorphic. It is trivial to see that there are continuously many
virtually non-isomorphic countable, discrete, abelian groups, for
instance by considering all groups $\Sigma_{n \in I} \Bbb Z/p_n
\Bbb Z_n$ with $I \subset \Bbb N$ and $p_n$ the prime numbers and
noticing that there are only countably many groups in each virtual
isomorphism class.
\hfill Q.E.D.

\proclaim{2.14. Corollary} Let $G$ be an infinite property
$({\text{\rm T}})$ group and $\sigma$ the Bernoulli shift action
of $G$ on $\Pi_g(\Bbb T, \lambda)_g$. Denote
$\sigma_n=\sigma^{\Bbb Z/n \Bbb Z}$, where $\Bbb Z/n\Bbb Z$ acts
as a (diagonal) product action on $\Pi_g(\Bbb T, \lambda)_g$ and
$\sigma^{\Bbb Z/n \Bbb Z}$ is defined out of $\sigma$ as in $2.7$.
Then $\sigma_n$ is not w-malleable, $\forall n\geq 2$. If in
addition $G$ is ICC then the inclusion of factors $N = A^{\Bbb Z/n
\Bbb Z} \rtimes_{\sigma_n} G \subset A \rtimes_\sigma G=M$ has
Jones index $[M:N] = n$ $\text{\rm ([J83])}$, with $M$ constructed from a
Bernoulli shift action of an ICC Kazhdan group, while $N$ cannot
be constructed from such data, i.e. $N$ cannot be realized as
$N=A_0 \rtimes_{\sigma_0} G_0$ with $G_0$  ICC Kazhdan group and
$(\sigma_0, G_0)$ a Bernoulli shift action.
\endproclaim
\vskip .05in {\it Proof}. If $\sigma_n$ would be w-malleable then
by Theorem 2.6 we would have $\text{\rm H}^1(\sigma_n,G) = {\text{\rm
Char}}(G)$. But by Theorem 2.9 we have $\text{\rm H}^1(\sigma_n,G)=\text{\rm
Char}(G) \times \Bbb Z/n\Bbb Z$ and since $G$ has (T), Char$(G)$
is finite so Char$(G) \ncong \text{\rm Char}(G) \times \Bbb Z/n\Bbb
Z$ for $n \geq 2$, a contradiction.

If $N = A_0 \rtimes_{\sigma_0} G_0$ for some Bernoulli shift
action $\sigma_0$ of an ICC property (T) group $G_0$ then by the
Superrigidity result (7.6 in [P04]) it would follow that
$(\sigma_n, G)$ and $(\sigma_0, G_0)$ are conjugate actions, with
$G \simeq G_0$, showing in particular that $\text{\rm
H}^1(\sigma_0, G_0) =\text{\rm H}^1(\sigma_n, G)$. Since
$\sigma_0$ is a $G_0$-Bernoulli shift action and $G_0\simeq G$ has
property (T), by ([PSa03]) we have $\text{\rm
H}^1(\sigma_0,G_0)={\text{\rm Char}}(G)$, while by Corollary 2.12
$\text{\rm H}^1(\sigma_n, G)= \text{\rm Char}(G) \times \Bbb
Z/n\Bbb Z$, a contradiction. \hfill Q.E.D.

\heading 3. 1-cohomology for actions of free product groups
\endheading

We now use Theorem 2.9 to calculate the 1-cohomology for quotients
of $G$-Bernoulli shifts in the case $G$ is a free product of
groups, $G=*_{n \geq 0} G_n$, with all $G_n$ either amenable or in
the class $w\Cal T$, at least one of them with this latter
property. Rather than locally compact as in 2.9,  the $\text{\rm H}^1$-group
is ``huge'' in this case, having either $\Cal U(A)$ or $
\Cal U(A)/\Bbb T$ as direct
summand. Since however $\Cal U(A), \Cal U(A)/\Bbb T$ are easily seen to be
connected, the quotient of $\text{\rm H}^1$ by the connected component of
$1$ provides a ``nicer'' group which is calculable and still an
invariant to stable orbit equivalence. This allows us to
distinguish many actions for each such $G$.

\proclaim{3.1. Lemma} Let $G_0, G_1, ...$ be a sequence of groups
and $G = *_{n \geq 0} G_n$ their free product.  Let $\sigma: G
\rightarrow \text{\rm Aut}(X, \mu)$ be a free ergodic measure
preserving action of $G$ on the probability space.

$1^\circ$. For each sequence $w=(w_i)_{i \geq 0}$ with $w_i \in
\text{\rm Z}^1(\sigma_{|G_i},G_i)$, $i \geq 0$, there exists a unique $\Delta(w)
\in \text{\rm Z}^1(\sigma,G)$ such that $\Delta(w)_{|G_i}=w_i, \forall i \geq
0$. The map $w \mapsto \Delta(w)$ is an isomorphism between the
Polish groups $\Pi_{i \geq 0} \text{\rm Z}^1(\sigma_{|G_i},G_i)$ and
$\text{\rm Z}^1(\sigma,G)$.

$2^\circ$. If $\sigma_{|G_0}$ is ergodic and
$K \subset \text{\rm Z}^1(\sigma_{|G_0})$ is a Polish subgroup which
maps $1$ to $1$ onto
its image $K/\sim_c$ in $\text{\rm H}^1(\sigma_{|G_0},G_0)$, then
$\Delta$
defined in $1^\circ$ implements a $1$ to $1$ continuous morphism $\Delta'$
from $K \times \Pi_{j \geq 1}
\text{\rm Z}^1(\sigma_{|G_j},G_j)$ into $\text{\rm H}^1(\sigma,G)$.
If in addition
$K/\sim_c=\text{\rm H}^1(\sigma_{|G_0},G_0)$ is also
surjective (so if the $1$-cohomology
exact sequence for $\text{\rm H}^1(\sigma_{|G_0},G_0)$ is split),
then $\Delta'$ is an
onto isomorphism of topological groups. In particular, if $\sigma_{|G_0}$
is weakly mixing then $\Delta$ implements a $1$ to $1$
continuous morphism $\Delta'$
from ${\text{\rm Char}} (G_0) \times \Pi_{j \geq 1}
\text{\rm Z}^1(\sigma_{|G_j},G_j)$ into $\text{\rm H}^1(\sigma,G)$ and if
$\text{\rm H}^1(\sigma_{|G_0},G_0)= {\text{\rm Char}}(G_0)$ then $\Delta'$ is an
onto isomorphism of Polish groups.
\endproclaim
\vskip .05in {\it Proof}. Part $1^\circ$ is evident by the
isomorphism between $\text{\rm Z}^1(\sigma,G)$ and Aut$_0(A \rtimes_\sigma
G;A)$ or by noticing that a function $w: G \rightarrow \Cal U(A)$
is a 1-cocycle for $\sigma$ iff $\{w_gu_g\}_{g \in G} \subset A
\rtimes_\sigma G$ is a representation of $G$.

$2^\circ$. If there would exist $u \in \Cal U(A)$ such that
$\Delta(w)_g = \sigma_g(u)u^*, \forall g,$ where $w = (w_0,
w_1, ...)$ for some $w_0 \in K$, $w_i \in
\text{\rm Z}^1(\sigma_{|G_i})$,$ i \geq 1$,
then $\sigma_{g_0}(u)u^* = w_0(g_0)$,
$\forall g_0 \in G_0$, implying that $w_0=1$.
By the ergodicity of $\sigma_{|G_0}$ this implies
$u \in \Bbb C1$. Thus $w_i=1, \forall i$, so that $w = (1,1,...,1)$.
If in addition
$\text{\rm H}^1(\sigma_{|G_0})\simeq K$ then $\Delta'$ follows
onto because $\Delta$ is onto and because of the way $\Delta$ is
defined. \hfill Q.E.D.

\vskip .05in

From 3.1.2$^\circ$ above we see that in case $\sigma_{|G_0}$ is
weakly mixing, then in order to calculate $\text{\rm H}^1(\sigma,
G)$ for $G=*_{n \geq 0} G_n$ we need to know $\text{\rm
H}^1(\sigma_{|G_0}, G_0)$ and $\text{\rm Z}^1(\sigma_{|G_i}, G_i)$
for $i\geq 1$. By 2.12, all these groups can be calculated if
$\sigma$ is the Bernoulli shift action of $G$, or certain
quotients of it, and $G_0\in w\Cal T$. The groups $\text{\rm Z}^1$
can in fact be calculated for amenable equivalence relations as
well, as shown below.

For convenience, we denote by $\Bbb G$ the
Polish group $\Cal U(A)$ (with the topology given by convergence
in norm $\|\cdot \|_2$), where $A=L^\infty(\Bbb T, \mu)$ as usual,
and by $\Bbb G_0$ the ``pointed'' space $\Bbb G/\Bbb T$. It is
easy to see that $\Bbb G$ is contractible (use for instance the
proofs in [PT93]), so that both $\Bbb G, \Bbb G_0$ are
connected. Also,  $\Bbb G^\infty \simeq \Bbb G$ and $\Bbb G
\times \Bbb G_0 \simeq \Bbb G_0$.

\proclaim{3.2. Lemma} $1^\circ$. If $G_i\in w\Cal T$,
$\Lambda$ is a countable discrete
abelian group and $\sigma'$ is an action of $G_i$ of the form
$\sigma_\Lambda$, as constructed in $2.10$, then $\text{\rm Z}^1(\sigma')
\simeq \Bbb G_0 \times \text{\rm Char}(G_i) \times \Lambda$.

$2^\circ$. If $\sigma'$ is a free m.p. action of a finite group
with $n \geq 2$ elements on the probability space $(X,\mu)$ and $Y
\subset X$ is a measurable subset with $\mu(Y)=(n-1)/n$ then
$\text{\rm H}^1(\sigma')=\{{\bold 1}\}$ and
$\text{\rm B}^1(\sigma')=\text{\rm Z}^1(\sigma')
\simeq \Cal U(L^\infty(Y,\mu))$. In particular, if $(X, \mu)$ is
non-atomic then $\text{\rm Z}^1(\sigma')\simeq \Bbb G$.

$3^\circ$. If $\sigma'$ is a free ergodic m.p.
action of an infinite amenable
group then $\text{\rm Z}^1(\sigma')\simeq \Bbb G$. Moreover,
$\text{\rm B}^1(\sigma')$ is
proper and dense in $\text{\rm Z}^1(\sigma')$.
\endproclaim
\vskip .05in \noindent {\it Proof}. Part $1^\circ$ is clear by
2.6, 2.9, 2.10 and the last part of 2.12, while $2^\circ$ is
folklore result.

$3^\circ$. By 1.4.1 and the results of Dye and Ornstein-Weiss
([D63], [OW80]), we may assume the infinite amenable group is
equal to $\Bbb Z$ and that the action is mixing (say a Bernoulli
shift). Identify $\text{\rm Z}^1(\sigma', \Bbb Z)$ with Aut$_0(A
\rtimes \Bbb Z, A)$ and notice that if $u=u_1\in M= A
\rtimes_{\sigma'} \Bbb Z$ denotes the canonical unitary
implementing the single automorphism $\sigma'(1)$ of $A$ then any
$v \in \Cal U(A)$ implements a unique automorphism $\theta^v \in
\text{\rm Aut}_0(M, A)$ satisfying $\theta^v(au)=avu$. Also, it is
trivial to see that $\Cal U(A) \ni v \mapsto \theta^v \in
\text{\rm Aut}_0(M, A)$ is an isomorphism of topological groups.
The fact that $\text{\rm B}^1(\sigma')$ is dense in $\text{\rm
Z}^1(\sigma')$ is immediate to deduce from ([OW80], [CFW81]) and
part $2^\circ$. Also, by 2.4.1$^\circ$ we have $\Bbb T \subset
\text{\rm H}^1(\sigma')$, so the subgroup $\text{\rm
B}^1(\sigma')$ is proper in $\text{\rm Z}^1(\sigma')$. \hfill
Q.E.D.

\proclaim{3.3. Theorem} Let $\{G_n\}_{n \geq 0}$ be a sequence
of groups, at least two of them non-trivial, and denote $G =
*_{n\geq 0} G_n$ their free product. Let $J=\{j \geq 0\mid G_j \in
w\Cal T \}$ and assume $0 \in J$ and $G_j$ amenable for all $j$
not in $J$. Let $\Lambda$ be a countable discrete abelian group
and denote by $\sigma_\Lambda$ the action of the group $G$
constructed in $2.12$, as a quotient of the classic $G$-Bernoulli
shift. We have the following isomorphisms of Polish groups:

$1^\circ$. If $J=\{0\}$ then $\text{\rm H}^1(\sigma_\Lambda, G)\simeq \Bbb G
\times \text{\rm Char}(G_0) \times \Lambda$.

$2^\circ$. If $J\neq \{0\}$ (i.e. $|J|\geq 2$) then
$\text{\rm H}^1(\sigma_\Lambda, G) \simeq \Bbb G_0^{|J|-1}
\times \Pi_{j \in J} \text{\rm
Char}(G_j) \times \Lambda^{|J|}$.
\endproclaim
\vskip .05in \noindent {\it Proof}. This is now trivial by 3.1,
3.2 and by the properties of $\Bbb G, \Bbb G_0$. \hfill Q.E.D.

\vskip .05in
\noindent
{\it 3.4. Notation}. Let $\sigma$ be a free
ergodic m.p. action of an infinite countable discrete group $G$ on a
standard probability space. We denote by $\tilde{\text{\rm
H}}^1(\sigma,G)$ the quotient of $\text{\rm Z}^1(\sigma,G)$ by the
connected component Z$^1_0(\sigma,G)$ of ${\bold 1}$ in
$\text{\rm Z}^1(\sigma,G)$. Since $\text{\rm
Z}^1_0(\sigma,G)$ is a closed subgroup in $\text{\rm
H}^1(\sigma,G)$, $\tilde{\text{\rm H}}^1(\sigma,G)$
with its quotient topology is a totally disconnected Polish group.
Note that, since B$^1(\sigma,G)$ is connected (being the immage of the
connected topological group $\Bbb G$), one has B$^1(\sigma,G) \subset
\text{\rm Z}^1_0(\sigma,G)$ and $\tilde{\text{\rm H}}^1(\sigma,G)$
coincides with the quotient of H$^1(\sigma,G)$ by the connected
component of ${\bold 1}$ in H$^1(\sigma,G)$.
Also, since $\text{\rm H}^1(\sigma,G)$ is invariant to stable
orbit equivalence, so is $\tilde{\text{\rm H}}^1(\sigma,G)$. If
$\Cal G$ is an ergodic full groupoid as in 1.3, then
$\tilde{H}^1(\Cal G)$ is defined similarly and has similar
properties.

\proclaim{3.5. Corollary} Under the same assumptions as in $3.3$,
if all $G_j, j \in J$, have finite character group (for instance
if they have the property $(\text{\rm T})$), or more generally if
$\text{\rm Char}(G_j)$ is totally disconnected $\forall j\in J$, then
$\tilde{\text{\rm H}}^1(\sigma_\Lambda, G) \simeq  \Pi_{j \in J} \text{\rm
Char}(G_j) \times \Lambda^{|J|}$ as Polish groups.
\endproclaim
\vskip .05in \noindent {\it Proof}. Trivial by 3.3 and the
comments in 3.4.
\hfill Q.E.D.

\proclaim{3.6. Corollary} Let $H_1, H_2, ..., H_k$ be infinite
property $(\text{\rm T})$ groups and $0 \leq n \leq \infty$. The
free product group $H_1 * H_2 * ... * H_k * \Bbb F_n$ has
uncountably many non stably orbit equivalent free ergodic m.p.
actions.
\endproclaim
\vskip .05in \noindent {\it Proof}. Clear by 3.5 and by the
argument in the proof of 2.13. \hfill Q.E.D.

\vskip .05in Note that the groups $G=*_{n \geq 0} G_n$ for which
we calculated the 1-cohomology for quotients of $G$-Bernoulli
shifts in this section do have infinite subgroups $H_0\subset G$
such that $(G,H_0)$ has the relative property (T): for instance,
if $H_0 \subset G_0$ is the infinite wq-normal subgroup of $G_0
\in w\Cal T$ such that $(G_0, H_0)$ has the relative property (T)
then $(G,H_0)$ has the relative property (T). It is trivial to see
though that $gG_0g^{-1} \cap G_0$ is finite $\forall g \in
G\setminus G_0$, so that by $(2.3')$ $H_0$ is not wq-normal in
$G$. Even more so, from 2.6 and 3.2, we deduce:

\proclaim{3.7. Corollary} If $G = K_1 * K_2$ with
$K_1, K_2$ non-trivial groups,
then $G$ is not in the class $w\Cal T$.
\endproclaim
\vskip .05in \noindent {\it Proof}. If $K_1, K_2$ are finite then
$G$ has the Haagerup approximation property ([H79]), so it cannot
contain an infinite subgroup with the relative property (T) (see
e.g. [P03]). If say $K_1$ is infinite and we let $\sigma$ be a
$G$-Bernoulli shift, then by 3.1.2$^\circ$ and $3.2$
H$^1(\sigma,G)$ contains either $\Bbb G$ or $\Bbb G_0$ as closed
subgroups. Since the latter are not compact (not even locally
compact), this contradicts 2.6. \hfill Q.E.D.

\vskip .05in \noindent {\bf 3.8. Remarks}. 1$^\circ$. Let
$\overline{\text{\rm H}}^1(\sigma,G)$ denote the quotient of
Z$^1(\sigma,G)$ by the closure of B$^1(\sigma,G)$ in
Z$^1(\sigma,G)$, or equivalently the quotient of H$^1(\sigma,G)$
by the closure of $\hat{\bold 1}$ in H$^1(\sigma,G)$. We see by
the definition that $\overline{\text{\rm H}}^1(\sigma,G)$ is
invariant to stable orbit equivalence. One can use arguments
similar to the ones in ([P01], [P03]) to prove that if $G\in w\Cal
T$ has an infinite amenable quotient $K$ with $\pi:G \rightarrow
K$ the quotient map, and $\sigma_g=\sigma_0(g) \otimes
\sigma_1(\pi(g))$, where $\sigma_0$ is a Bernoulli shift action of
$G$ and $\sigma_1$ a Bernoulli shift action of $K$, then
$\overline{\text{\rm H}}^1(\sigma,G)=\text{\rm Char}(G)$, while
$\sigma$ is not strongly ergodic in this case. By using the
construction in the proof of Theorem 2.9,  from the action
$\sigma$ one can then construct free ergodic m.p. actions
$\sigma_\Lambda$ of $G$ such that $\overline{\text{\rm
H}}^1(\sigma_\Lambda,G)=\text{\rm Char}(G)\times \Lambda$, for any
countable abelian groups $\Lambda$.

2$^\circ$. Corollary 2.12 and Theorem 3.3 provide computations of
the 1-cohomology group H$^1(\sigma_\Lambda, G)$ for the family of
actions $\sigma_\Lambda$ constructed in 2.12, for most groups $G$
having infinite subgroups with the relative property (T). However,
groups having the Haagerup compact approximation property ([H79]),
such as the free groups $\Bbb F_n, 2 \leq n \leq \infty$, do not
contain infinite subgroups with the relative property (T) (see
e.g. [P02]). The problem of calculating the H$^1$-groups for
$G$-Bernoulli shifts and their quotients $\sigma_\Lambda$ when $G$
are free groups, or other non-amenable groups with the Haagerup
property, remains open. Note however that by 3.1.1$^\circ$ and
3.2.3$^\circ$ if $\sigma$ is an arbitrary free ergodic m.p. action
of $\Bbb F_n$ on the probability space then Z$^1(\sigma, \Bbb
F_n)\simeq \Cal U(A)^n=\Bbb G^n \simeq \Bbb G$, so
$\tilde{\text{\rm H}}^1(\sigma,G)=\{1\}$. (In fact Z$^1$ is even
contractible.) Also, by 3.1.2$^\circ$ one has an embedding of
$\Bbb T  \times \Bbb G^{n-1}$ into H$^1(\sigma, \Bbb F_n)$
whenever one of the generators of $\Bbb F_n$ acts weak mixing
(e.g. when $\sigma$ is a Bernoulli shift). All this indicates that
the H$^1$-invariant may be less effective in recognizing orbit
inequivalent actions of the free groups.

Related to this, our last result below emphasizes the limitations
of the ``deformation/rigidity'' techniques of ([P01], [P03]) when
trying to prove that $\text{\rm H}^1(\sigma,G)=\text{\rm Char}(G)$
for commutative and non-commutative $G$-Bernoulli shifts, beyond
the class of w-rigid groups $G$ dealt with in ([P01], [PoS03]) and
the class $w\Cal T$ in this paper.

Thus, we let this time $G$ be an arbitrary non-amenable group and
$\sigma$ be the action of $G$ on the finite von Neumann algebra
$(N,\tau)= \overline{\otimes}_g (N_0,\tau_0)_g$ by (left)
Bernoulli shift automorphisms, with the ``base'' $(N_0,\tau_0)$
either the diffuse abelian von Neumann algebra $L^\infty(\Bbb
T,\lambda)$, or a finite dimensional factor $M_{n \times n}(\Bbb
C)$, or the hyperfinite $\text{\rm II}_1$ factor $R$. By ([P01],
[P03]) $\sigma$ is malleable. More precisely, there exists a
continuous action $\alpha$ of $\Bbb R$ on $(N\overline{\otimes} N,
\tau\otimes \tau)$ such that $[\alpha,\tilde{\sigma}]=0$ and
$\alpha_1(N \otimes 1) = 1\otimes N$, where $\tilde{\sigma}_g =
\sigma_g \otimes \sigma_g, g\in G$.

\proclaim{3.8. Proposition} Let $w \in \text{\rm Z}^1(\sigma,G)$.
The following conditions are equivalent:

$(i)$. $w$ is cohomologous to a character of $G$.

$(ii)$. For sufficiently small $|t|$, the representation $\pi_t$
of $G$ on $L^2(N, \tau)\overline{\otimes} L^2(N, \tau)$
given by $f \mapsto (w_g \otimes
1)\tilde{\sigma}_g(f)\alpha_t(w_g^*\otimes 1)$ is a direct sum
between a multiple of the trivial representation of $G$ and a
subrepresentation  of a multiple of the left regular
representation of $G$.
\endproclaim
\vskip .05in {\it Proof}. If $w_g =\sigma_g(u)u^* \gamma(g), g\in
G,$ for some $\gamma \in {\text{\rm Char}}(G)$ and $u \in \Cal
U(N)$, then $U_t(f) = (u\otimes 1) f \alpha_t(u^*\otimes 1)$, for
$f\in L^2(N,\tau)\overline{\otimes} L^2(N,\tau)$ defines a unitary
operator that intertwines the representations $\pi_0$ and $\pi_t$,
which thus follow equivalent. But $\pi_0=\tilde{\sigma}$ is a
direct sum between one copy of the trivial representation of $G$
and a subrepresentation of a multiple of the left regular
representation of $G$ (see e.g. [S80], [J83b]). This shows that
$(i) \implies (ii)$.

Conversely, since $\underset t  \rightarrow 0 \to \lim
\|\pi_t(g)(1) - 1\|_2=0, \forall g,$ where $1=1_N \otimes 1_N$, if
$\pi_t$ satisfy $(ii)$, then for $t$ small enough $\pi_t(g)(1)$
follows close to $1$ uniformly in $g \in G$, i.e. $(w_g \otimes 1)
\alpha_t(w_g^* \otimes 1), g \in G,$ is uniformly close to $1$ (in
the norm $\|\cdot \|_2$). But then the argument in ([P01]) or
([PoS03]) shows that $w$ is coboundary, thus $(ii) \implies (i)$.
\hfill Q.E.D.

\head  References\endhead

\item{[B91]} M. Burger, {\it Kazhdan constants for} $SL(3,\Bbb
Z)$, J. reine angew. Math., {\bf 413} (1991), 36-67.

\item{[C75]} A. Connes: {\it Outer conjugacy classes of
automorphisms of factors}, Ann. \'Ec. Norm. Sup. {\bf 8} (1975),
383-419.

\item{[C75b]} A. Connes: {\it Sur la classification des facteurs
de type} II, C. R. Acad. Sci. Paris {\bf 281} (1975), 13-15.

\item{[C80]} A. Connes: {\it A type II$_1$ factor with countable
fundamental group}, J. Operator Theory {\bf 4} (1980), 151-153.

\item{[CFW81]} A. Connes, J. Feldman, B. Weiss: {\it An amenable
equivalence relation is generated by a single transformation},
Ergodic Theory Dynamical Systems {\bf 1} (1981), 431-450.

\item{[CW80]} A. Connes, B. Weiss: {\it Property} (T) {\it and
asymptotically invariant sequences}, Israel J. Math. {\bf 37}
(1980), 209-210.

\item{[CoZ]} M. Cowling, R. Zimmer: {\it Actions of lattices in}
$Sp(n,1)$, Ergod. Th. Dynam. Sys. {\bf 9} (1989), 221-237.

\item{[D63]} H. Dye: {\it On groups of measure preserving
transformations} II, Amer. J. Math, {\bf 85} (1963), 551-576.

\item{[FM77]} J. Feldman, C.C. Moore: {\it Ergodic equivalence
relations, cohomology, and von Neumann algebras I, II}, Trans.
Amer. Math. Soc. {\bf 234} (1977), 289-324, 325-359.

\item{[Fu99]} A. Furman: {\it Orbit equivalence rigidity}, Ann. of
Math. {\bf 150} (1999), 1083-1108.

\item{[G00]} D. Gaboriau: {\it Cout des r\'elations
d'\'equivalence et des groupes}, Invent. Math. {\bf 139} (2000),
41-98.

\item{[G02]} D. Gaboriau: {\it Invariants $\ell^2$ de r\'elations
d'\'equivalence et de groupes},  Publ. Math. I.H.\'E.S. {\bf 95}
(2002), 93-150.

\item{[GP03]} D. Gaboriau, S. Popa: {\it An Uncountable Family of
Non Orbit Equivalent Actions of $\Bbb F_n$}, math.GR/0306011, to
appear in Journal of the AMS.

\item{[Ge87]} S.L. Gefter: {\it On cohomologies of ergodic actions
of a T-group on homogeneous spaces of a compact Lie group}
(Russian), in ``Operators in functional spaces and questions of
function theory'', Collect. Sci. Works, Kiev, 1987, pp 77-83.

\item{[GeGo88]} S.L. Gefter, V.Y. Golodets: {\it Fundamental
groups for ergodic actions and actions with unit fundamental
groups}, Publ RIMS {\bf 6} (1988), 821-847.

\item{[dHV89]} P. de la Harpe, A. Valette: ``La propri\'et\'e T de
Kazhdan pour les groupes localement compacts'', Ast\'erisque {\bf
175}, Soc. Math. de France (1989).

\item{[H79]} U. Haagerup: {\it An example of a non-nuclear}
$C^*$-{\it algebra which has the metric approximation property},
Invent. Math. {\bf 50} (1979), 279-293.

\item{[Hj02]} G. Hjorth: {\it A converse to Dye's Theorem}, UCLA
preprint, September 2002.

\item{[Jol02]} P. Jolissaint: {\it On the relative property T}, to
appear in l'Ens. Math.

\item{[J83]} V.F.R. Jones : {\it Index for subfactors}, Invent.
Math. {\bf 72} (1983), 1-25.

\item{[J83b]} V.F.R. Jones : {\it A converse to Ocneanu's
theorem}, J. Operator Theory, {\bf 10} (1983), 121-123.

\item{[JS87]} V.F.R. Jones, K. Schmidt: {\it Asymptotically
invariant sequences and approximate finiteness}, Amer. J. Math
{\bf 109} (1987), 91-114.

\item{[K67]} D. Kazhdan: {\it Connection of the dual space of a
group with the structure of its closed subgroups}, Funct. Anal.
and its Appl. {\bf 1} (1967), 63-65.

\item{[L94]} A. Lubotzky: ``Discrete groups, expanding graphs and
invariant measures'', Birkh\"{a}user, 1994.

\item{[Ma82]} G. Margulis: {\it Finitely-additive invariant
measures on Euclidian spaces}, Ergodic. Th. and Dynam. Sys. {\bf
2} (1982), 383-396.

\item{[MoSh02]} N. Monod, Y. Shalom: {\it Orbit equivalence
rigidity and bounded cohomology}, Preprint 2002.

\item{[M82]} C.C. Moore: {\it Ergodic theory and von Neumann
algebras}, Proc. Symp. Pure Math. {\bf 38} (Amer. Math. Soc.
1982), pp. 179-226.

\item{[MvN43]} F. Murray, J. von Neumann: {\it Rings of operators
IV}, Ann. Math. {\bf 44} (1943), 716-808.

\item{[P01]} S. Popa: {\it Some rigidity results for
non-commutative Bernoulli shifts}, MSRI preprint 2001-005.

\item{[P02]} S. Popa: {\it On a class of type II$_1$ factors with
Betti numbers invariants}, MSRI preprint no 2001-024, revised
math.OA/0209310, to appear in Ann. of Math.

\item{[P03]} S. Popa: {\it Strong Rigidity of} II$_1$ {\it Factors
Arising from Malleable Actions of $w$-Rigid Groups} I,
math.OA/0305306.

\item{[P04]} S. Popa: {\it Strong Rigidity of} II$_1$ {\it Factors
Arising from Malleable Actions of $w$-Rigid Groups} II,
math.OA/0407137.

\item{[PSa03]} S. Popa, R. Sasyk: {\it On the cohomology of
actions of groups by Bernoulli shifts}, math.OA/0310211.

\item{[PT93]} S. Popa, M. Takesaki: {\it The topological structure
of the unitary and automorphism groups of a factor},
Communications in Math. Phys. {\bf 155} (1993), 93-101.

\item{[S80]} K. Schmidt: {\it Asymptotically invariant sequences
and an action of $SL(2, \Bbb Z)$ on the $2$-sphere}, Israel. J.
Math. {\bf 37} (1980), 193-208.

\item{[S81]} K. Schmidt: {\it Amenabilty, Kazhdan's property T,
strong ergodicity and invariant means for ergodic group-actions},
Ergod. Th. \& Dynam. Sys. {\bf 1} (1981), 223-236.

\item{[Si55]} I.M. Singer: {\it Automorphisms of finite factors},
Amer. J. Math. {\bf 177} (1955), 117-133.

\item{[Z84]} R. Zimmer: ``Ergodic Theory and Semisimple Groups'',
Birkhauser, Boston, 1984.

\enddocument